\documentclass[11pt,a4paper,reqno]{amsproc}

\usepackage{amsfonts}
\usepackage{amsmath,amscd}
\usepackage{amsthm}

{ \theoremstyle{plain}
\newtheorem{theorem}{Theorem}
\newtheorem{lemma}{Lemma}
\newtheorem{corollary}{Corollary}
\newtheorem{proposition}{Proposition}
}

{ \theoremstyle{definition}
\newtheorem{definition}{{\bf Definition}}}

{ \theoremstyle{definition}

}

{ \theoremstyle{remark}
\newtheorem{remark}{Remark}
}

\DeclareMathOperator{\ord}{ord}

\DeclareMathOperator{\supp}{\mathrm{supp}}
\DeclareMathOperator{\wind}{\mathrm{wind}}
\DeclareMathOperator{\Mod}{\mathrm{mod}}
\DeclareMathOperator{\Res}{\mathcal{R}}
\DeclareMathOperator*{\res}{\mathrm{res}}

\def\Com#1{\mathbb{C}^{#1}}

\def\I{\mathrm{i}}

\def\scal#1#2{\langle #1, #2\rangle }

\def\Log{{\rm Log\,}}
\def\D{{\mathbb D}}

\def\C{{\mathbb C}}
\def\Z{{\mathbb Z}}

\def\Wilde{}

\title[The resultant on compact Riemann surfaces]
{The resultant on compact Riemann surfaces
 }

\author{Bj\"orn Gustafsson}

\address{Mathematical department, KTH}
\email{gbjorn@kth.se}

\author{Vladimir G. Tkachev}
\address{Mathematical Department,\\Volgograd State University\\\textsl{Current:} Mathematical department, KTH}
\email{tkatchev@kth.se}



\begin{document}

\begin{abstract}
We introduce a notion of resultant of two
meromorphic functions on a compact Riemann surface and demonstrate its usefulness
in several respects. For example, we exhibit several integral formulas for
the resultant, relate it to potential theory
and give explicit formulas for the algebraic dependence between two
meromorphic functions on a compact Riemann surface. As a particular application,
the exponential transform of a quadrature domain in the
complex plane is expressed in terms of the resultant of two
meromorphic functions on the Schottky double of the domain.
\end{abstract}

\date{}

\maketitle


\section{Introduction}

A bounded domain $\Omega$ in the complex plane
is called a (classical) {\it quadrature domain}  \cite{Aharonov-Shapiro76},
\cite{Sakai82}, \cite{Shapiro92}, \cite{Gustafsson-Shapiro05}
or, in a different terminology, an {\it algebraic domain} \cite{Varchenko-Etingof94},
if there exist finitely many points $z_i\in \Omega$ and coefficients $c_{kj}\in \C$
($i=1,\dots, N$, say) such that
\begin{equation}\label{QD1}
\int_\Omega h\, dxdy =\sum_{k=1}^N \sum_{j=1}^{s_k} c_{kj} h^{(j-1)}(z_k)
\end{equation}
for every integrable analytic function $h$ in $\Omega$ \cite{qd}.
In the last two decades there has been a growing interest in the
applications of quadrature domains to various problems in mathematics
and theoretical physics, ranging from Laplacian growth to integrable
systems and string theory (see recent articles
\cite{Gustafsson-Putinar07}), \cite{KKMWZ}, and the references therein).

One of the most intriguing properties of quadrature domains is their
algebraicity \cite{Aharonov-Shapiro76}, \cite{Gustafsson83}: the
boundary of a quadrature domain is (modulo finitely many points) the
full real section of an algebraic curve:
\begin{equation}\label{lamma}
\partial \Omega = \{z\in \mathbb{C}: Q(z,\bar z)=0\},
\end{equation}
where $Q(z,w)$ is an irreducible Hermitian polynomial. Moreover, the
corresponding full algebraic curve (essentially $\{(z,w)\in
\mathbb{C}^2: Q(z,w)=0\}$) can be naturally identified with the
Schottky double $\widehat{\Omega}$ of $\Omega$ by means of the
Schwarz function $S(z)$ of $\partial\Omega$. The latter satisfies
$S(z)=\bar{z}$ on $\partial\Omega$ and is, in the case of a
quadrature domain, meromorphic in all $\Omega$.
%

A deep impact into the theory of quadrature domains was the
discovery  by M.~Putinar \cite{Putinar96} in the mid 1990's of an
alternative characterization in terms of hyponormal operators.
Recall that J.~Pincus proved \cite{Pincus} that with any bounded
linear operator $T:H\to H$ in a Hilbert space $H$ for which the
self-commutator is positive (i.e., $T$ is hyponormal) and has rank
one, say
\begin{equation*}
[T^*,T]=T^*T-TT^*=\xi\otimes \xi,  \qquad 0\ne \xi\in H,
\end{equation*}
one can associate a unitary invariant, the so-called principal
function. This is  a measurable function $g:\Com{}\to [0,1]$,
supported on the spectrum of $T$, such that for any $z,w$ in the
resolvent set of $T$ there holds
\begin{equation}\label{pincus11}
\det (T_z^*T_w{T_z^*}^{-1}{T_w}{}^{-1})=\exp[\frac{1}{2\pi\I}\int_\C
\frac{g(\zeta)\;d\zeta \wedge d\bar{\zeta}}{(\zeta-z)(\bar{\zeta}
-\bar {w})}],
\end{equation}
where $T_u=T-uI$. The right hand side in (\ref{pincus11}) is
referred to as the \textit{exponential transform} of the function
$g$. In case $g$ is the characteristic function of a bounded set
$\Omega$ we have the exponential transform of $\Omega$,
\begin{equation}\label{emem}
E_\Omega (z,w) =\exp[\frac{1}{2\pi\I}\int_\Omega \frac{d\zeta}{\zeta
-z}\wedge \frac{d\bar{\zeta}}{\bar{\zeta} -\bar {w}}].
\end{equation}
A central result in Putinar's theory is the following criterion: a
domain $\Omega$ is a quadrature domain if and only if the
exponential transform of $\Omega$ is a rational function of the form
\begin{equation}\label{e111}
E_\Omega(z,w)=\frac{Q(z, w)}{P(z)\overline{P(w)}}, \qquad |z|, |w|\gg 1,
\end{equation}
where $P$ and $Q$ are polynomials. In this case $Q$ is the same as the polynomial in (\ref{lamma}).


In the present paper we shall unify the above pictures by
interpreting the exponential transform of a quadrature domain in
terms of resultants of meromorphic functions on the Schottky double
of the domain. To this end we need to extend the classical concept
of resultant of two polynomials to a notion of resultant for
meromorphic functions on a compact Riemann surface. The introduction
of such a meromorphic resultant and the demonstration of its
usefulness in several contexts is the main overall purpose of this
paper.

The definition of the resultant is natural and simple: given two
meromorphic functions $f$ and $g$ on a compact Riemann surface $M$
we define their \textit{meromorphic resultant} as
\begin{equation*}
\Res(f,g)=\prod_{i=1}^m\frac{g(a_i)}{g(b_i)},
\end{equation*}
where $(f)=\sum a_i-\sum b_i=f^{-1}(0)-f^{-1}(\infty)$ is the
divisor of $f$. This resultant actually depends only on the divisors
of $f$ and $g$.
It follows from Weil's reciprocity law that the resultant is symmetric:
\begin{equation*}
\Res(f,g)=\Res(g,f).
\end{equation*}

For the genus zero case the meromorphic resultant is just a
cross-ratio product of four polynomial resultants, whereas  for
higher genus surfaces it can be expressed as a cross-ratio product
of values of theta functions. In the other direction, the classical
resultant of two polynomials (which can be viewed as meromorphic
functions with a marked pole) may be recovered from the meromorphic
one by specifying a local symbol at the infinity (see
Section~\ref{sec:local}).

It is advantageous in many contexts to amplify the resultant to an
\textit{elimination function}. With $f$ and $g$ as above this is
defined as
\begin{equation*}
\mathcal{E}_{f,g}(z,w)=\Res(f-z, g-w),
\end{equation*}
where $z,w$ are free complex parameters.
Thus defined, $\mathcal{E}_{f,g}(z,w)$ is a \textit{rational}
function in $z$ and $w$
having the elimination property
\begin{equation*}
\mathcal{E}_{f,g}(f(\zeta),g(\zeta))= 0 \qquad (\zeta\in M).
\end{equation*}
In particular, this gives an explicit formula for the algebraic
dependence between two meromorphic functions on a compact Riemann
surface. Treating the variables $z$ and $w$ in the definition of
$\mathcal{E}_{f,g}(z,w)$ as spectral parameters in the elimination
problem, the above identity resembles the Cayley-Hamilton theorem
for the characteristic polynomial in linear algebra. This analogy
becomes more clear by passing to the so-called differential
resultant in connection with the spectral curves for two commutating
ODE's, see for example \cite{Pr91}.


The above aspects of the resultant and the elimination function
characterize them essentially from an algebraic side. In the paper
we shall however much emphasize the analytic point of view by
relating the resultant to objects such as the exponential transform
(\ref{emem}) and the Fredholm determinant. One of the key results is
an integral representation of the resultant
(Theorem~\ref{integralformula}), somewhat similar to (\ref{emem}).
From this we deduce one of the main results of the paper: the
exponential transform of a quadrature domain $\Omega$ coincides with
a natural elimination function on the Schottky double
$\widehat{\Omega}$ of $\Omega$:
\begin{equation}\label{later}
E_\Omega(z,w)= \mathcal{E}_{f,f^*}(z,\bar w).
\end{equation}
Here $(f, f^*)$ is a canonical pair of meromorphic functions on
$\widehat{\Omega}$: $f$  equals the identity function on $\Omega$,
which extends to a meromorphic function on the double
$\widehat{\Omega}$ by means of the Schwarz function, and $f^*$ is
the conjugate of the reflection of $f$ with respect to the
involution on $\widehat{\Omega}$. In Section~\ref{seq:appl1} we use
formula (\ref{later}) to construct explicit examples of classical
quadrature domains.

In Section~\ref{sec:sze} we discuss the meromorphic resultant
$\Res(f,g)$ as a function of the quotient
\begin{equation}\label{polar}
h(z)=\frac{f(z)}{g(z)}.
\end{equation}
Clearly, $f$ and $g$ are not uniquely determined by $h$ in this
representation, but given $h$ it is easy to see that there are, up
to constant factors, only finitely many pairs $(f,g)$ with non-zero
resultant $\Res(f,g)$ for which (\ref{polar}) holds.
Thus, the natural problem of characterizing the total range
$\sigma(h)$ of these values $\Res(f,g)$ arises.

Another case of interest is that the divisors of $f$ and $g$ are
confined to lie in prescribed \textit{disjoint} sets. This makes
$\Res(f,g)$ uniquely determined by $h$ and connects the subject to
classical work of E.~Bezout and L.~Kronecker on representations of
the classical resultant $\Res_{\mathrm{pol}}(f,g)$ by
Toeplitz-structured determinants with entries equal to Laurent
coefficients of the quotient $h(z)$. The 60's and 70's  brought
renewed interest to this area in connection with asymptotic behavior
of truncated Toeplitz determinants for rational generating functions
(cf. \cite{Baxter}, \cite{day}, \cite{FH}). This problem naturally
occurs in statistical mechanics in the study of the spin–spin
correlations for the two-dimensional Ising model (see, e.g.,
\cite{Botch99}) and  in quantum many body systems \cite{Forr},
\cite{Bas}.

One of the general results for rational symbols is an exact formula
given by M.~Day \cite{day} in 1975. Suppose that $h$ is a rational
function with simple zeros which is regular on the unit circle and
does not vanish at the origin and infinity: $\ord_0 h\leq 0$,
$\ord_\infty h\leq 0$. Then for any $N\geq 1$:
\begin{equation*}
\det (h_{i-j})_{1\leq i,j\leq N}=\sum_{i=1}^{p} r_i H_i^N,
\qquad h_k=\frac{1}{2\pi }
\int_{0}^{2\pi}e^{-\I k\theta}h(e^{\I \theta})d\theta,
\end{equation*}
where $p$, $r_i$, $H_i$ are suitable rational expressions
in the divisor of $h$.
An accurate analysis of these expressions  reveals the
following interpretation
of the above identity in terms of resultants:
\begin{equation}\label{Beeth}
\frac{\det (h_{i-j})_{1\leq i,j\leq N}}{h^N(0)}=\sum \Res(z^Nf,g),
\end{equation}
where the (finite) sum is taken over all pairs $(f,g)$ satisfying
(\ref{polar}) such that $g$ is normalized by $g(\infty)=1$ and the
divisor of zeros of $g$ coincides with the restriction of the polar
divisor of $h$ to the unit disk: $ (g)_+=(h)_-\cap \mathbb{D}$.

In the above notation, the equality (\ref{Beeth}) can be thought of
as an identity between the elements of $\sigma(h)$ with a prescribed
partitioning of the divisor. In Section~\ref{seq:residen} we
consider resultant identities in the genus zero case in general, and
show that there is a family of linear relations on $\sigma(h)$.
These identities may be formally interpreted as a limiting case (for
$N=0$) of the above Day formula (\ref{Beeth}). Moreover, our
resultant identities are similar to those given recently by
A.~Lascoux and P.~Pragacz \cite{LascPr} for Sylvester's double sums.
On the other hand, by specializing the divisor $h$ we obtain a
family of trigonometric identities generalizing known trigonometric
addition theorems. Some of these identities were obtained recently
by F.~Calogero in \cite{Calog,Calog1}. For non-zero genus surfaces
the situation with describing $\sigma(h)$ becomes much more
complicated. We consider some examples for a complex torus, which
indicates a general tight connection between resultant identities
and addition theorems for theta-functions.

Returning to (\ref{later}) and comparing this identity  with
determinantal representation (\ref{pincus11}) we find it reasonable
to conjecture that one can associate to any compact Riemann surface
an appropriate functional calculus for which the elimination
function becomes a Fredholm determinant. In Section~\ref{sec:seg} we
demonstrate such a model for the zero genus case. We show that the
meromorphic resultant of two rational functions is given by a
determinant of a multiplicative commutator of two Toeplitz operators
on an appropriate Hardy space. There are interesting similarities
between our determinantal representation (cf. formula (\ref{both})
below) of the meromorphic resultant and the tau-function for
solutions of some integrable hierarchies (see, for instance,
\cite{SEG}).

Further aspects of the meromorphic resultant discussed in the paper
are interpretations in terms of potential theory, in
Section~\ref{sec5}, and various cohomological points of view, e.g.,
an expression of the resultant in terms of the Serre duality pairing
(subsections~\ref{linebundle} and ~\ref{subsec:serre}). In
Section~\ref{seq:int} we give an independent proof of the symmetry
of the resultant using the formalism of currents, and also derive
several integral representations. Section~\ref{sec3} contains the
main definitions and other preliminary material, and in
Section~\ref{sec:pol} we review the polynomial resultant.

The authors are grateful to Mihai Putinar, Emma Previato and Yurii Neretin
for many helpful
comments and to the Swedish Research Council and the Swedish Royal Academy of
Sciences for financial support.  This research is a part of the European Science
Foundation Networking Programme ``Harmonic and Complex Analsyis and
Applications HCAA''.

\section{The polynomial resultant}
\label{sec:pol}

The resultant  of two  polynomials, $f$ and $g$,
in one complex variable is a polynomial function in the coefficients of $f$, $g$
having the elimination property that it vanishes if and only if $f$ and $g$ have
a common zero \cite{Waerden}. The resultant is a  classical concept which goes back
to the work of L.~Euler, E.~B\'{e}zout, J.~Sylvester and A.~Cayley. Traditionally, it plays an important role in algorithmic
algebraic geometry as an effective tool for elimination of variables in polynomial equations. The renaissance of the classical
theory of elimination in the last decade owes much to recent progress in toric geometry, complexity theory and the theory of univariate and  multivariate
residues of rational forms (see, for instance, \cite{Gelfand-Kapranov-Zelevinskij}, \cite{St97}, \cite{tsih}, \cite{Dickenstein}).

We begin with some basic definitions and facts. In terms of the zeros
of polynomials
\begin{equation}\label{PQ}
f(z)=f_m\prod_{i=1}^m(z-a_i)=\sum_{i=0}^m f_iz^i ,
\quad g(z)=g_n\prod_{j=1}^n(z-c_j) =\sum_{j=0}^n g_jz^j,
\end{equation}
the resultant is given by the Poisson product formula
\cite{Gelfand-Kapranov-Zelevinskij}
\begin{equation}\label{res1}
\begin{split}
\Res_{\mathrm{pol}}(f,g)&=f_m^ng_n^m\prod_{i,j} (a_i-c_j) =f_m^n\prod_{i=1}^m g(a_i)=(-1)^{mn}g_n^m\prod_{j=1}^n f(c_j).
\end{split}
\end{equation}
It follows immediately
from this definition that $\Res_{\mathrm{pol}}(f,g)$  is skew-sym\-met\-ric
and multiplicative:
\begin{equation}\label{mulmul}
\Res_{\mathrm{pol}}(f,g)=(-1)^{mn}\Res_{\mathrm{pol}}(g,f), \qquad \Res_{\mathrm{pol}}(f_1f_2,g)=\Res_{\mathrm{pol}}(f_1,g)\Res_{\mathrm{pol}}(f_2,g).
\end{equation}

Alternatively, the resultant is uniquely (up to a normalization)
defined as the irreducible integral polynomial in the coefficients
of $f$ and $g$ which vanishes if and only if $f$ and $g$ have a
common zero.


All known explicit representations of the  polynomial resultant
appear as certain determinants in the coefficients of the polynomials. Below we briefly comment on the most important
determinantal representations. The interested reader may consult the recent monograph
\cite{Gelfand-Kapranov-Zelevinskij} and the surveys \cite{Dickenstein},
\cite{St97}, where further information on the
subject can be found.

With $f$, $g$ as above, let us define an operator
$
S:\mathcal{P}_{n}\oplus \mathcal{P}_{m}\to \mathcal{P}_{m+n}
$
by the rule:
\begin{equation*}\label{sysl}
S(X,Y)=fX+gY,
\end{equation*}
where $\mathcal{P}_{k}$ denotes the space of polynomials of degree
$\leq k-1$ ($\dim \mathcal{P}_k=k$). Then
\begin{equation}\label{sfor}
\Res_{\mathrm{pol}}(f,g)=\det
\begin{pmatrix}
  f_0   &        &       & g_0          &        &&&\\
  \vdots  &\ddots&       & \vdots          & \ddots &&&\\
  f_m   &        & f_0   & g_n          &        &&&g_0    \\
        & \ddots & \vdots&              &\ddots  &&& \vdots   \\
        &        & f_m   &              &        &&& g_n\\
\end{pmatrix}
\end{equation}
where the latter is the \textit{Sylvester matrix} representing
$S$ with respect to the monomial basis.


An alternative method to describe the resultant is the so-called
B\'{e}zout-Cayley formula.
For $\deg f=\deg g=n$ it reads
\begin{equation*}\label{bezou}
\Res_{\mathrm{pol}}(f,g)=\det (\beta_{ij})_{0\leq i,j\leq n-1},
\end{equation*}
where
\begin{equation}\label{bezout}
\frac{f(z)g(w)-f(w)g(z)}{z-w}=\sum_{i,j=0}^{n-1}\beta_{ij}z^{i}w^{j},
\end{equation}
is the B\'{e}zoutian of $f$ and $g$.
The general case, say $\deg f<\deg g$, is obtained from
(\ref{mulmul}) and (\ref{bezout}) by completing $f(z)$ to $z^kg(z)$,
$k=\deg g-\deg f$.

Other remarkable  representations of the resultant are given as
determinants of Toeplitz-structured matrices with entries equal to
Laurent coefficients of the quotient $h(z)=\frac{f(z)}{g(z)}$. These
formulas were known already to E.~Bezout and were rediscovered and
essentially developed  later by J.~Sylvester and L.~Kronecker in
connection to finding of the greatest common divisor of two polynomials
(see Chapter~12 in \cite{Gelfand-Kapranov-Zelevinskij} and \cite{Utesh1}).

Recently, a similar formula in terms of contour integrals  of the quotient $h(z)$ has been given by R.~Hartwig \cite{Hart} (see also M.~Fisher and R.~Hartwig  \cite{FH}). In its simplest form this formula reads as follows. With $f$ and $g$ as in (\ref{PQ}), we assume $g_0=g(0)\ne 0$.
Then for any $N\geq n$, the polynomial resultant, up to a constant factor, is the truncated
Toeplitz determinant for the symbol  $h(z)$:
\begin{equation}\label{hankel}
\begin{split}
\Res_{\mathrm{pol}}(f,g)
&=f_m^{n-N}g_0^{m+N}\det t_{m,N}(h),\\
\end{split}
\end{equation}
where $h(z)=\sum_{k=0}^\infty h_k z^k$ is the Taylor development of the quotient around $z=0$ and
\begin{equation*}
t_{m,N}(h)= \left(
              \begin{array}{cccc}
                h_{m} & h_{m-1} & \ldots & h_{m-N+1} \\
                h_{m+1} & h_{m} & \ldots & h_{m-N+2} \\
                \vdots & \vdots & \ddots & \vdots \\
                h_{m+N-1} & h_{m+N-2} & \ldots & h_{m} \\
              \end{array}
            \right),
\end{equation*}
and $h_k=0$ for negative $k$.


The determinant $\det t_{m,N}(h)$ is a commonly used object in
theory of Toeplitz operators. For instance, the celebrated
Szeg\"o limit theorem (see, e.g., \cite{Botch99}) states that,
under some natural assumptions, $\det t_{0,N}(h)$ behaves like a
geometric progression. Exact formulations will be given in
Section~\ref{seq:seg}, where the above identity is generalized
to the meromorphic case.

It is worth mentioning here another powerful and rather unexpected
application of $\det t_{m,N}(h)$, the so-called  Thom-Porteous
formula in the theory of determinantal varieties \cite{FultonP},
\cite[p.~415]{Griffith-Harris}. We briefly describe this identity
in the classical setup. Consider an $n\times m$ ($n\leq m$) matrix
$A$ with entries $a_{ij}$ being homogeneous forms in the variables
$x_1,\ldots, x_k$  of degree $p_i+q_j$ (for some integers $p_i$, $q_j$).
Denote by $V_r$ the locus of points in $\mathbb{P}^k$ at which the rank
of $A$ is at most $r$. Then, thinking of $p_i$, $q_j$ as formal variables,
one has
\begin{equation*}\label{porteo}
\deg V_r=\det t_{m-r,n-r}(c),\qquad \sum_{k=0}^\infty c_kz^k
=\frac{\prod_{j=1}^m(1+q_jz)}{\prod_{i=1}^n(1-p_iz)}.
\end{equation*}

We mention here also a differential analog of the polynomial
resultant in algebraic theory of commuting (linear) ordinary
differential operators. A key observation goes back to
J.L.~Burchnall and  T.W.~Chaundy and states that commuting ordinary
differential operators satisfy an equation for a certain algebraic
curve, the so-called spectral curve of the corresponding operators
(see  \cite{Pr96} for a detailed discussion and historical remarks).
The defining equation of the curve is equivalent to the vanishing of
a determinant of a Sylvester-type matrix. This phenomenon was a main
ingredient of the modern fundamental algebro-geometric approach
initiated by I.~Krichever \cite{Kr} in the theory of integrable
equations. By using the Burchnall-Chaundy-Krichever correspondence
between meromorphic functions on a suitable Riemann surface and
differential operators, E.~Previato in \cite{Pr91} succeeded to get
a pure algebraic version of the proof of Weil's reciprocity.

All the determinantal formulas given above fit into a general scheme:
given a pair of polynomials one can associate an operator $S$ in
a suitable coefficient model space such that
$\Res_{\mathrm{pol}}(f,g)=\det {S}.$ On the other hand, none of the
models behaves well under multiplication of polynomials. This makes
it difficult to translate identities like (\ref{mulmul}) into matrix
language. One way  to get around this difficulty is to observe that
(\ref{hankel}) is a special case of the Szeg\"o strong limit theorem
for rational symbols \cite{FH} and to consider infinite dimensional
determinantal (Fredholm) models instead. We sketch such a model in
Section~\ref{sec:seg} below.


\section{The meromorphic resultant}
\label{sec3}

\subsection{Preliminary remarks}
\label{sec:prel}
For rational functions with neither zeros nor poles at infinity, say
\begin{equation}\label{fgdefnew}
f(z)=\lambda\prod_{i=1}^m\frac{z-a_i}{z-b_i},
\quad g(z)=\mu\prod_{j=1}^n\frac{z-c_j}{z-d_j},
\end{equation}
($\lambda,\mu\ne 0$ and all $a_i,b_i,c_j,d_j$ distinct) it is
natural to define the resultant as
\begin{equation}
\label{gagb1new}
\Res(f,g)=\prod_{i=1}^m \frac{g(a_i)}{g(b_i)}
=\prod_{j=1}^n\frac{f(c_j)}{f(d_j)}.
\end{equation}
In other words,
\begin{equation}\label{crossr}
\begin{split}
\Res(f,g)&=\prod_{i=1}^m \prod_{j=1}^n\frac{a_i-c_j}{a_i-d_j}\cdot \frac{b_i-d_j}{b_i-c_j}
=\prod_{i=1}^m \prod_{j=1}^n (a_i,b_i,c_j,d_j),\\
\end{split}
\end{equation}
where $(a,b,c,d):=\frac{a-c}{a-d}\cdot \frac{b-d}{b-c}$
is the classical cross ratio of four points.

Note that (nonconstant) polynomials do not fit into this picture
since they always have a pole at infinity, but the polynomial
resultant can still be recovered by a localization procedure
(see Section~\ref{sec:local}). Notice also that the above resultant
for rational functions actually has better properties than the
polynomial resultant, e.g., it is symmetric ($\Res(f,g)=\Res(g,f)$),
homogenous of degree zero and it only depends on the divisors of
$f$ and $g$. The resultant for meromorphic functions on a compact
Riemann surface will be modeled on the above definition
\eqref{gagb1new} and contain it as a special case.

\subsection{Divisors and their actions}
We start with a brief discussion of divisors. A divisor on a
Riemann surface $M$ is a finite formal linear combination
of points on $M$, i.e., an expression of the form
\begin{equation}
\label{DD}
D=\sum_{i=1}^m n_i a_i,
\end{equation}
$a_i\in M$, $n_i\in \Z$. Thus a divisor is the same thing as a
$0$-chain,  which acts on $0$-forms, i.e., functions, by
integration. Namely, the divisor \eqref{DD} acts on functions
$\varphi$ by
\begin{equation}
\label{Dvarphi}
\scal{D}{\varphi}=\int_D \varphi=\sum_{i=1}^m n_i \varphi(a_i).
\end{equation}

From another (dual) point of view divisors can be looked upon as maps $M\to \Z$
with support at a finite number of points, namely the maps which evaluate the
coefficients in expressions like \eqref{DD}. If $D$ is a divisor as in \eqref{DD}
we also write $\Wilde{D}: M\to \Z$ for the corresponding evaluation map. Then $D=\sum_{a\in M} \Wilde{D}(a)a.$
The degree of $D$ is
\begin{equation*}
 \deg D= \sum_{i=1}^m n_i=\sum_{a\in M} \Wilde{D}(a).
\end{equation*}
and its support is
\begin{equation*}
\supp D=\{a\in M: D(a)\ne 0\}.
\end{equation*}

If $f:M\to \mathbb{P}$ is a nonconstant meromorphic function and $\alpha \in \mathbb{P}$ then
the inverse image $f^{-1}(\alpha)$, with multiplicities counted, can be considered
as a (positive) divisor in a natural way. The divisor of $f$ then is
\begin{equation}\label{princf}
(f)= f^{-1}(0)-f^{-1}(\infty).
\end{equation}
If $f$ is constant, not 0 or $\infty$, then $(f)=0$ (the zero element in the Abelian group of divisors).

Recall that any divisor of the form (\ref{princf}) is called a
\textit{principal} divisor. In the dual picture the same divisor
acts on points as follows:
\begin{equation*}
\Wilde{(f)} (a)=\ord_a (f),
\end{equation*}
where $\ord_a (f)$ is the integer $m$ such that, in terms of a local coordinate $z$,
\begin{equation*}
f(z)=c_m(z-a)^m + c_{m+1}(z-a)^{m+1}+\dots \quad {\rm with} \quad c_m\ne 0.
\end{equation*}
By $\ord f$ we denote the order of $f$, that is the cardinality of $f^{-1}(0)$.

Divisors act on functions by \eqref{Dvarphi}. We can also let
functions act on divisors. In this case we shall, by convention, let
the action be multiplicative rather than additive: if
$h=h(u_1,\ldots,u_k)$ is a function and $D_1,\ldots,D_k$ are divisors,
we set
\begin{equation}\label{hhh}
h(D_1,\dots, D_k)= \prod_{a_1,\ldots,a_k\in M} h(a_1,\dots,a_k)^{\Wilde{D}_1(a_1)\cdots \Wilde{D}_k(a_k)},
\end{equation}
whenever this is well-defined. Observe that this definition is
consistent with the standard evaluation of a function at a point.
Indeed, any point $a\in M$ may be regarded simultaneously as a
divisor $D_a=a$. Then $h(a_1,\ldots,a_p)=h(D_{a_1},\ldots,D_{a_p})$.
In what follows we make no distinction between $D_a$ and $a$.

With branches of the logarithm chosen arbitrarily (\ref{hhh}) can also be written
\begin{equation*}
h(D_1,\ldots,D_p)=
\exp \;\scal{D_1\otimes \ldots \otimes D_p}{\log h}.
\end{equation*}
When $D_i$, $i=1,\ldots,p$ are principal divisors, say $D_i=(g_i)$
for some meromorphic functions $g_i$, the definition (\ref{hhh}) yields
\begin{equation*}
h((g_1),\ldots,(g_p))=\prod_{a_1,\ldots,a_p\in M}h(a_1,\ldots,a_p)^{\ord_{a_1}(g_1)\cdots \ord_{a_p}(g_p)}.
\end{equation*}

\subsection{Main definitions}
Let now $f$, $g$ be meromorphic functions (not identically $0$ and $\infty$) on an
arbitrary compact Riemann surface $M$ and let their divisors be
\begin{equation}\label{gdiv}
\begin{split}
(f)& =f^{-1}(0)-f^{-1}(\infty)=\sum\nolimits_{i=1}^m a_i -\sum\nolimits_{i=1}^m b_i,\\
(g)&=g^{-1}(0)-g^{-1}(\infty)=\sum\nolimits_{j=1}^n c_j -\sum\nolimits_{j=1}^n d_j.
\end{split}
\end{equation}

At first we assume that $(f)$ and $(g)$ are ``generic'' in the sense of having disjoint supports.
In view of the suggested
resultant \eqref{gagb1new} for rational functions the following definition is
natural.

\begin{definition}
The (meromorphic) \textit{resultant} of two generic meromorphic functions $f$ and $g$ as above is
\begin{equation}\label{main1new}
\Res(f,g)=g((f))=\prod_{i=1}^m\frac{g(a_i)}{g(b_i)}
=\frac{g(f^{-1}(0))}{g(f^{-1}(\infty))}= \exp\scal{(f)}{ \log g}.
\end{equation}
In the last expression, an arbitrary branch of $\log g$ can be chosen at each point of $(f)$.
\end{definition}

Elementary properties of the resultant are multiplicativity in each variable:
\begin{equation*}
\Res(f_1 f_2, g)=\Res(f_1,g)\Res(f_2,g), \qquad \Res(f,g_1 g_2)=\Res(f,g_1)\Res(f,g_2).
\end{equation*}
An important observation is homogeneity of degree zero
\begin{equation}\label{homm}
\Res(af,bg)=\Res(f,g)
\end{equation}
for $a,b\in \Com{*}:=\C\setminus \{0\}$. The latter implies that $\Res(f,g)$ depends merely on the divisors $(f)$ and $(g)$.

Less elementary, but still true, is the symmetry:
\begin{equation}
\label{symmWnew}
\Res(f,g)=\Res(g,f),
\end{equation}
i.e., in the terms of the divisors
\begin{equation*}
\prod_i \frac{g(a_i)}{g(b_i)}=\prod_j \frac{f(c_j)}{f(d_j)}.
\end{equation*}
This is a consequence of Weil's reciprocity law \cite{Weil}, \cite[p.~242]{Griffith-Harris}.
In Section~\ref{seq:int} we shall find some integral formulas for the resultant and also give
an independent proof of  \eqref{symmWnew}.

If, in (\ref{hhh}), some of the divisors $D_k$ are principal then
the resulting action $h$ may be written as a composition of the
corresponding  resultants. For instance, for a function $h$ of two
variables we have
\begin{equation}\label{subseq}
h((f),(g))=\Res_u(f(u),\Res_v(g(v),h(u,v))),
\end{equation}
where $\Res_u$ denotes the resultant in the $u$-variable.

\begin{remark}
The definition of meromorphic resultant naturally extends to more general objects
than meromorphic functions. Indeed, of $f$ we need only its divisor
and $g$ may be  a fairly arbitrary function. We shall still use (\ref{main1new}) as a definition in such extended contexts. However, there is no symmetry relation like (\ref{symmWnew}) in general. See e.g. Lemma~\ref{rem:conj}.

\end{remark}

When, as above, $(f)$ and $(g)$ have disjoint supports $\Res(f,g)$
is a nonzero complex number. It is important to extend the
definition of $\Res(f,g)$ to certain cases when $(f)$ and $(g)$ do
have common points.

\begin{definition}\label{defadm}
A pair of two meromorphic functions $f$ and $g$ is said to be
\textit{admissible} on a set $A\subset M$ if the function $a\to
\ord_a(g)\ord_a(f)$ is sign semi-definite on $A$
(i.e., is either $\geq 0$ on all $A$ or $\leq 0$ on all $A$). If
$A=M$ we shall simply say that $f$ and $g$ is an admissible pair.
\end{definition}

It is easily seen that the product in (\ref{main1new}) is
well-defined as a complex number or $\infty$ whenever $f$ and $g$ form an admissible pair.

Clearly, any pair of two meromorphic functions whose divisors have
no common points is admissible  (we call such pairs generic).
Another  important example  is  the
family of all polynomials, regarded as meromorphic
functions on the Riemann sphere $\mathbb{P}$. It is easily seen that
any pair of polynomials is admissible with respect to an
\textit{arbitrary} subset $A\subset\mathbb{P}$.

The following elimination property is an immediate corollary of the definitions.

\begin{proposition}\label{pr:cr}
Let two nonconstant meromorphic functions $f$, $g$ form an admissible pair on $M$. Then $\Res (f,g)=0$ if and only if $f$ and $g$ have  a common zero or a common pole.
In particular,
$
\Res(f,g)=0$
if $f$ and $g$ are polynomials.

\end{proposition}

\subsection{Elimination function}
We have seen above that the meromorphic resultant of two \textit{individual} functions is not always
well-defined (namely, if the two functions do not form an admissible pair). However one may still get useful information by embedding the functions in families depending on parameters, for example by taking the resultant of $f-z$ and $g-w$. We shall see in Section~\ref{uvidim} that such resolved
versions of the resultant have additional analytic advantages.

Let $z,w\in \C$ be free variables. The  expression
\begin{equation*}
\mathcal{E}(z,w)\equiv \mathcal{E}_{f,g}(z,w)=\Res(f-z, g-w),
\end{equation*}
if defined, will be called the \textit{elimination function} of $f$ and $g$.

\begin{theorem}\label{th:elim}
Let $f$ and $g$ be nonconstant meromorphic functions without common
poles. Then the elimination function is well defined everywhere
except for finitely many pairs $(z,w)$, and it is  a rational function
of the form
\begin{equation*}\label{ee1}
\mathcal{E}(z,w)=\frac{Q(z,w)}{P(z)R(w)},
\end{equation*}
where $Q$, $P$, $R$ are polynomials, and
\begin{equation*}
P(z)= \prod_{d\in g^{-1}(\infty)} (z-f(d)), \qquad
R(w) = \prod_{b\in f^{-1}(\infty)} (w-g(b)).
\end{equation*}
\end{theorem}

\begin{proof}
Note that a linear transformation $f\to f-z$ keeps the polar locus
unchanged. Thus the elimination function $\Res(f-z, g-w)$ is
well-defined for all pairs $(z,w)$ such that $f^{-1}(z)\cap
g^{-1}(\infty)=g^{-1}(w)\cap f^{-1}(\infty)=\emptyset$. Let $(z,w)$
be any such pair. Then applying the symmetry relation
(\ref{symmWnew}) we obtain
\begin{equation*}
\mathcal{E}(z,w)=\frac{(g-w)(f^{-1}(z))}{(g-w)(f^{-1}(\infty))}
= \frac{(f-z)(g^{-1}(w))}{(f-z)(g^{-1}(\infty))}.
\end{equation*}

Let $f$, $g$ have orders $m$ and $n$, respectively, as in (\ref{gdiv}), and let $\{f_i^{-1}\}$ denote the branches
of $f^{-1}$. Then spelling out the meaning we find, using that the symmetric functions of
$\{g(f_i^{-1}(z))\}$ are single-valued from the Riemann sphere into itself,
hence are rational functions, that
\begin{equation*}
(g-w)(f^{-1}(z)) = \prod_{i=1}^m (g(f_i^{-1}(z))-w)=(-1)^m (w^m + R_1(z) w^{m-1}+\dots + R_m(z)),
\end{equation*}
where the $R_i(z)$ are rational. Similarly,
\begin{equation*}
(g-w)(f^{-1}(\infty))
=(-1)^m (w^m + r_1 w^{m-1}+\dots + r_m),
\end{equation*}
where the $r_i$ are constants.

With the same kind of arguments for $(f-z)(g^{-1}(w))$ and
$(f-z)(g^{-1}(\infty))$ we obtain
\begin{equation*}
\begin{split}
\mathcal{E}(z,w)&=\frac{w^m + R_1(z) w^{m-1}+\dots + R_m(z)}{w^m + r_1 w^{m-1}+\dots + r_m}
=\frac{z^n + P_1(w) z^{n-1}+\dots + P_n(w)}{z^n + p_1 z^{n-1}+\dots + p_n}.
\end{split}
\end{equation*}
Clearing the denominators (in the numerators) yields the required statement.
\qed\end{proof}

Important, and useful in applications, is the following elimination
property of the function $\mathcal{E}_{f,g}(z,w)$. Let us choose
$\zeta\in M$ arbitrarily and insert $z=f(\zeta)$, $w=g(\zeta)$ into
$\mathcal{E}_{f,g}(z,w)$. Since the functions $f-z$ and $g-w$ then have a common
zero (namely at $\zeta$) this gives, by Proposition~\ref{pr:cr}, that
\begin{equation*}
\mathcal{E}_{f,g}(f(\zeta), g(\zeta))=0 \quad (\zeta\in M).
\end{equation*}
In particular,
\begin{equation*}
Q(f,g)=0,
\end{equation*}
i.e., we have recovered the classical polynomial relation
between two functions on a compact Riemann surface (see
\cite{Farkas-Kra}, \cite{Forster}, for example).

\subsection{Extended elimination function}
We have seen  that the elimination function is well-defined for any
pair of meromorphic functions without common poles. One step
further, linear fractional transformations allow us to refine the
definition of elimination function in such a way that it becomes
well-defined for \textit{all} pairs of meromorphic functions.

Namely, let $f$ and $g$ be two arbitrary meromorphic  functions and consider the function of four complex variables:
\begin{equation}\label{utvidg}
\mathcal{E}(z,w;z_0,w_0)\equiv \mathcal{E}_{f,g}(z,w;z_0,w_0)=\Res\biggl(\frac{f-z}{f-z_0}, \frac{g-w}{g-w_0}\biggr).
\end{equation}

Let us choose  arbitrary the pair $(z,z_0)$. Then we have for divisor: $
(\frac{f-z}{f-z_0})=f^{-1}(z)-f^{-1}(z_0)$. It is easy to see that the resultant in (\ref{utvidg}) is well defined for any quadruple  $(z,w;z_0,w_0)$ with
\begin{equation}\label{holds1}
[g^{-1}(w)\cup g^{-1}(w_0)]\cap [f^{-1}(z)\cup f^{-1}(z_0)]=\emptyset.
\end{equation}
The set $X$ of all $(z,w;z_0,w_0)$ such that (\ref{holds1}) holds is a dense open subset of in $\Com{4}$.

Applying then an argument similar to that in Theorem~\ref{th:elim}, we find
that the right hand side in (\ref{utvidg}) is a rational function for $(z,w;z_0,w_0)\in X$. We call this function the \textit{extended} elimination function of $f$ and $g$.

We have the cross-ratio-like symmetries
$\mathcal{E}(z,w;z_0,w_0)=\mathcal{E}(z_0,w_0;z,w)$, and
\begin{equation*}\label{sym2}
\mathcal{E}(z,w_0;z_0,w)=\frac{1}{\mathcal{E}(z,w;z_0,w_0)}.
\end{equation*}

In the case when the elimination function $\mathcal{E}_{f,g}(z,w)$ is well-defined we have the following reduction:
\begin{equation*}\label{quadr-e}
\mathcal{E}(z,w;z_0,w_0)=\frac{\mathcal{E}(z,w)\mathcal{E}(z_0,w_0)}{\mathcal{E}(z,w_0)\mathcal{E}(z_0,w)}=\frac{Q(z,w)Q(z_0,w_0)}{Q(z,w_0)Q(z_0,w)},
\end{equation*}
with $Q$ as in  Theorem~\ref{th:elim}.

In the other direction, the ordinary elimination function, if well-defined,
can be viewed as a limiting case of the extended version. Indeed,
it follows from null-homogeneity of the meromorphic resultant that
\begin{equation*}
\mathcal{E}(z,w;z_0,w_0)=\Res\biggl(\frac{f-z}{1-f/z_0}, \frac{g-w}{1-g/w_0}\biggr),
\end{equation*}
and therefore that
\begin{equation*}
\lim_{z_0,w_0\to\infty}\mathcal{E}(z,w;z_0,w_0)=\mathcal{E}(z,w).
\end{equation*}

There are still cases when the elimination function is not defined or is trivial while its extended version contains information. To illustrate this, let us consider a meromorphic function $f$ of order $n$ and let $g=f$. Then a straightforward computation reveals that
\begin{equation*}
\mathcal{E}_{f,f}(z,w;z_0,w_0)=\left(\frac{z-z_0}{z-w_0}\cdot \frac{w-w_0}{w-z_0}\right)^n=(z,w,z_0,w_0)^n,
\end{equation*}
where $(z,w,z_0,w_0)$ is the cross ratio.

\subsection{The meromorphic resultant on surfaces with small genera}
\label{seq:ratt}  

On the Riemann sphere $\mathbb{P}$ the resultant reduces to a product of cross ratios
(\ref{crossr}) and the symmetry relation (\ref{symmWnew}) becomes
trivial. Note that the cross ratio itself may be regarded as
the meromorphic resultant of two linear fractional functions.

From a computational point of view, evaluation of the meromorphic resultant on $\mathbb{P}$ is similar to the evaluation of polynomial resultants. Indeed, for any admissible rational functions given by the ratio of polynomials,  $f=f_1/f_2$ and $g=g_1/g_2$, one finds that
\begin{equation}\label{def1}
\Res(f,g)=f(\infty)^{\ord_\infty(g)}g(\infty)^{\ord_\infty(f)}\cdot\frac{\Res_{\mathrm{pol}}(f_1,g_1)\Res_{\mathrm{pol}}(f_2,g_2)}
{\Res_{\mathrm{pol}}(f_1,g_2)\Res_{\mathrm{pol}}(f_2,g_1)}.
\end{equation}
The latter formula combined with formulas in Section~\ref{sec:pol} expresses the meromorphic
resultant in terms of the \textit{coefficients} of the representing polynomials of $f$ and $g$.
For example, since each resultant in (\ref{def1}) is a Sylvester determinant (\ref{sfor}),
$$
\Res_{\mathrm{pol}}(f_i,g_j)=\det S(f_i,g_j)\equiv \det S_{ij},
$$
the resulting product amounts to
\begin{equation*}\label{matrx}
\Res(f,g)=f(\infty)^{\ord_\infty(g)}g(\infty)^{\ord_\infty(f)}\cdot\det (S_{12}^{-1}S_{11}S_{21}^{-1}S_{22}).
\end{equation*}

In Section~\ref{sec:seg} we give another, more invariant,
approach to the representation of meromorphic resultants via
determinants (see also Section~\ref{sec:Cauchy} for the
exponential representations of $\Res(f,g)$).


Now we spell out the definition of the resultant in case of Riemann surfaces of genus one. Consider the complex torus $M=\Com{}/L_\tau$, where $L_\tau=\mathbb{Z}+\tau \mathbb{Z}$ is the lattice formed by $\tau\in \Com{}$, $\mathrm{Im}\,\tau >0$. A meromorphic function on $M$ is represented as an $L_\tau$-periodic function on $\Com{}$. Let
\begin{equation*}
\theta(\zeta)=\theta_{11}(\zeta)\equiv \sum_{k=-\infty}^\infty e^{\pi \I (k^2\tau +k(1+\tau+2\zeta))}
\end{equation*}
be the Jacobi theta-function. Then  any meromorphic function $f$ on $M$ is given by a ratio of translated theta-functions:
\begin{equation*}\label{jacobi}
f(\zeta)=\lambda \prod_{i=1}^m\frac{ \theta (\zeta-a_i)}{\theta (\zeta-b_i)},
\end{equation*}
and a necessary and sufficient condition that such a ratio really defines a meromorphic function is that the divisor is principal, i.e., by Abel's theorem, that
\begin{equation}\label{jacobcond}
\sum_{i=1}^m (a_i-b_i)\in L.
\end{equation}

With $f$ as above and $g$ similarly with $c_j$ and $d_j$, $\sum_{j=1}^n(c_j-d_j)\in L$,  the following representation for the meromorphic resultant on the torus holds:
\begin{equation*}\label{jacobi1}
\Res(f,g)=\prod_{i=1}^m\prod_{j=1}^n\frac{ \theta (c_j-a_i)\theta (d_j-b_i)}{\theta (c_j-b_i)\theta (d_j-a_i)}.
\end{equation*}


\section{Integral representations}\label{seq:int}

\subsection{Integral formulas}
\label{sec:begin}
We shall derive some integral representations for the meromorphic resultant,
and in passing also give a proof of the symmetry \eqref{symmWnew},
Weil's reciprocity law. Let $f$, $g$ be nonconstant meromorphic functions on a
compact Riemann surface $M$ of genus $p\geq 0$ and recall \eqref{main1new} that the
resultant can be written
\begin{equation*}
\Res(f,g)= \exp\scal{(f)}{ \log g}.
\end{equation*}
We assume that the divisors $(f)$ and $(g)$ have disjoint supports.
Since $(f)$ is integer-valued and different branches of $\log g$
differ by integer multiples of $2\pi \I$ it does not matter which
branch of $\log g$ is chosen at each point of $(f)$. However, our
present aim is to treat $\log g$ as a global object on $M$, in order
to interpret $\scal{(f)}{\log g}$ as a current acting on a function
and to write it as an integral over $M$.

First of all, to any divisor $D$ can be naturally associated a
$2$-form current $\mu_D$ (a $2$-form with distribution
coefficients), which represents $D$ in the sense that
\begin{equation*}
\scal{D}{\varphi} =\int_D\varphi = \int_M\varphi \wedge \mu_D
\end{equation*}
for smooth functions $\varphi$. With $D=\sum n_i {a_i}$ this $\mu_D$
is of course just
\begin{equation}\label{mud}
\mu_D =\delta_D dx\wedge dy =\sum n_i \delta_{a_i} dx\wedge dy,
\end{equation}
where $\delta_a$ is the Dirac delta at the point $a$ and with
respect to a local variable $z=x+\I y$ chosen (only $\delta_a
dx\wedge dy$ has an invariant meaning). When $D=(f)$ we have the
following formula.

\begin{lemma}\label{omegaf}
If $f$ is a meromorphic function, then
$
\mu_{(f)}=\frac{1}{2\pi \I}\;d(\frac{df}{f})
$
in the sense of currents.
\end{lemma}

\begin{proof}
In a neighbourhood of a point $a$ with $\ord_a (f)=m$, i.e.,
\begin{equation*}
f(z)=c_m(z-a)^m + c_{m+1}(z-a)^{m+1}+\dots, \quad c_m\ne 0,
\end{equation*}
in terms of a local coordinate, we have
$
\frac{df}{f}= (\frac{m}{z-a}+h(z))dz
$
with $h$ holomorphic. Hence,
\begin{equation*}
d\biggl(\frac{df}{f}\biggr)=\frac{\partial}{\partial \bar{z}}\left(\frac{m}{z-a}+h(z)\right)\;d\bar{z}\wedge dz
=m\pi \delta_a d\bar{z}\wedge dz = 2\pi \I m\delta_a dx\wedge dy,
\end{equation*}
from which the lemma follows.
\qed
\end{proof}

Next we shall make $\log f$ and $\log g$ single-valued on $M$ by
making ``cuts''.
Let  $\alpha_1$,\dots,  $\alpha_p$, $\beta_1$,\dots,  $\beta_p$ be
a canonical homology basis for $M$ such that each $\beta_k$
intersects $\alpha_k$ once from the right to the left ($k=1,\dots, p$)
and no other crossings occur. We may choose these curves so that they do not
meet the divisors $(f)$ and $(g)$.

Since the divisors $(f)$ and $(g)$ have degree zero we can write
\begin{equation*}
(f)=\partial \gamma_f, \quad (g)=\partial \gamma_g
\end{equation*}
where $\gamma_f$, $\gamma_g$ are $1$-chains. We may arrange these
curves so that there are no intersections  and so that they are contained in $M\setminus(\alpha_1
\cup\dots\cup \beta_p)$.

Now, it is possible to select single-valued branches
of $\log f$ and  $\log g$ in
\begin{equation*}
M^\prime = M\setminus( \gamma_f \cup\gamma_g \cup \alpha_1 \cup\dots\cup \beta_p).
\end{equation*}
Fix such branches and denote them $\Log f$, $\Log g$.
Then $\Log f$ and $\Log g$ are functions, defined almost everywhere
on $M$, and $\Log g$ is smooth in a neighbourhood of the support of
$(f)$ and vice versa. In particular, $\scal{(f)}{\Log g}$ and $\scal{(g)}{\Log f}$ make sense.

Now using Lemma~\ref{omegaf} and partial
integration (with exterior derivatives taken in the sense of
currents) we get
\begin{equation*}
\begin{split}
\Res(f,g) &= \exp\scal{(f)}{\Log g} =\exp[\int_{M} \mu_{(f)} \wedge \Log g]\\
&=\exp[\frac{1}{2\pi \I}\int_{M} d(\frac{df}{f}) \wedge \Log g]
=\exp[\frac{1}{2\pi \I}\int_M \frac{df}{f}\wedge d \, {\Log \,} g].
\end{split}
\end{equation*}

In summary:
\begin{theorem}\label{integralformula}
Let $f$ and $g$ be two meromorphic functions on a compact Riemann
surface whose divisors have disjoint supports. Then
\begin{equation*}
\Res (f,g)
= \exp [\frac{1}{2\pi \I}\int_M \frac{df}{f}\wedge d \, {\Log \,} g].
\end{equation*}
In particular, for generic $z,w$,
\begin{equation*}
\mathcal{E}_{f,g}(z,w)= \exp[\frac{1}{2\pi \I}\int_M \frac{df}{f-z}\wedge d\, \Log (g-w)].
\end{equation*}
\end{theorem}

It should be noted that the only contributions to the integrals above
come from the jumps of $\Log g$ (and $\Log (g-w)$ respectively), because outside this set of
discontinuities the integrand contains $dz\wedge dz=0$ as a factor.

\subsection{Symmetry of the resultant}
We proceed to study $d\,\Log $ in detail. Let first $a,b$ be two points in the complex plane and
$\gamma$ a curve from $b$ to $a$ such that $\partial \gamma= a-b$
(formal difference).
Then, with a single-valued branch of the logarithm chosen in
$\C\setminus \gamma$,
\begin{equation*}
\begin{split}
d\,\Log \frac{z-a}{z-b} &=\frac{dz}{z-a}-\frac{dz}{z-b}
+\I [d\,\text{Arg\,}\frac{z-a}{z-b}]_{\text{jump contribution from $\gamma$}}
\\
&=\frac{dz}{z-a}-\frac{dz}{z-b} -2\pi \I dH_\gamma (z).
\end{split}
\end{equation*}
Here $dH_\gamma$ is the $1$-form current supported by $\gamma$ and
defined as the (distributional) differential of the function
$H_\gamma$ which in a neighbourhood of any interior point of
$\gamma$ equals $+1$ to the right of $\gamma$ and zero to the left.
Thus $dH_\gamma$ is locally exact  away from the end points.
The function $H_\gamma$ cannot be defined in any full neighbourhood
of $a$ or $b$. On the other hand, $dH_\gamma$ is taken to have no
distributional contributions at $a$ and $b$. One easily checks that
this gives a current which represents $\gamma$ in the sense that
\begin{equation*}
\int_\gamma \tau=\int_M dH_\gamma \wedge\tau
\end{equation*}
for all smooth $1$-forms $\tau$.
Taking $\tau$ of the form $d\varphi$ gives
\begin{equation*}
\int_M d(dH_\gamma)\wedge \varphi
=\int_M dH_\gamma\wedge d\varphi =\int_\gamma d\varphi
=\int_{\partial\gamma} \varphi.
\end{equation*}
Thus the $0$-chain, or divisor,
$\partial\gamma$ is represented by
$d(dH_\gamma)$. We can write this also as
$ 
d(dH_\gamma)=\mu_{\partial \gamma},
$ 
where $\mu_D$ is defined in (\ref{mud}).
Note in particular that $dH_\gamma$ is not closed, despite the notation.

If $\gamma$ and $\sigma$ are two curves (1-chains) which cross each other
at a point $c$, then it is easy to check (and well-known)
that
\begin{equation*}
dH_\gamma\wedge dH_\sigma =\pm \delta_c \;dx\wedge dy,
\end{equation*}
with the plus sign if $\sigma$ crosses $\gamma$ from the right (of $\gamma$)
to the left, the minus sign in the opposite case.
For the curves $\alpha_1,\dots,\beta_p$ in the canonical homology basis,
the forms $dH_{\alpha_1},\dots , dH_{\beta_p}$  are closed, since the curves are themselves closed.

Now we extend the above analysis to $\Log f$ in place of $\Log \frac{z-a}{z-b}$.
In addition to the jump across $\gamma_f$ (an arbitrary 1-chain in $M\setminus (\alpha_1\cup\ldots\cup \beta_p)$ with $\partial\gamma_f=(f)$) we need to take into account possible jumps
across the $\alpha_k$, $\beta_k$. In order to reach the right hand side of
$\alpha_k$ from the left hand side within $M^\prime$ one just follows $\beta_k$.
The increase of $\Log f$ along this curve is $\int_{\beta_k} \frac{df}{f}$, hence
this is also the jump of $\Log f$ across $\alpha_k$, from the left to the right.
With a similar analysis for the jump across $\beta_k$ one arrives at the following
expression for $d\,\Log f$:
\begin{equation*}
d\,\Log f = \frac{df}{f} - 2\pi \I (dH_{\gamma_f} +\sum_{k=1}^p
(\frac{1}{2\pi \I }\int_{\beta_k} \frac{df}{f} \cdot dH_{\alpha_k}
-\frac{1}{2\pi \I }\int_{\alpha_k} \frac{df}{f} \cdot
dH_{\beta_k})).
\end{equation*}
This means that $\gamma_f$ needs to be modified to the 1-chain
\begin{equation*}\label{sigmaf}
\sigma_f =\gamma_f  +\sum_{k=1}^p( \wind_{\beta_k}(f) \cdot{\alpha_k}
-\wind_{\alpha_k}(f)\cdot{\beta_k}),
\end{equation*}
where, for a closed curve $\alpha$ in general, $\wind_{\alpha}(f)$ stands for the winding number
\begin{equation*}
\wind_{\alpha}(f)= \frac{1}{2\pi \I }\int_\alpha \frac{df}{f} \in\Z.
\end{equation*}
Notice that $\partial\sigma_f=\partial\gamma_f=(f)$ and that now $\Log f$ can be taken to be single-valued analytic in $M\setminus \supp \sigma_f$.
The above can be we can summarized as follows.

\begin{lemma}\label{lem2}
Given any meromorphic function $f$ in $M$ there exists a 1-chain $\sigma_f$ having the property that $\partial \sigma_f=(f)$, $\log f$ has a single-valued branch, $\Log f$, in $M\setminus \supp \sigma_f$ and the exterior differential of $\Log f$, regarded as a 0-current in $M$ with jumps taken into account, is
\begin{equation*}
d\,\Log f=\frac{df}{f} - 2\pi \I dH_{\sigma_f}.
\end{equation*}
\end{lemma}

Since $\frac{df}{f}\wedge\frac{dg}{g}=0$ the lemma combined with Theorem~\ref{integralformula} gives the
following alternative formula for the resultant.

\begin{corollary} With notations as above
\begin{equation}\label{res2}
\Res(f,g) = \exp(-\int_M \frac{df}{f}\wedge dH_{\sigma_g})=
\exp \int_{\sigma_g} \frac{df}{f}.
\end{equation}
\end{corollary}

In the corollary $\sigma_f$ may be replaced by any 1-chain $\gamma$ with $\partial \gamma=(g)$, because this will make a difference in the integral only by an integer multiple of $2\pi\I$.

Next we compute
\begin{equation*}
\begin{split}
d\,\Log f &\wedge d\,\Log g
=(\frac{df}{f} - 2\pi \I dH_{\sigma_f} )\wedge (\frac{dg}{g} -2\pi \I dH_{\sigma_g} )
\\
=& \frac{df}{f}\wedge d\, \Log g
+ d\,\Log f \wedge \frac{dg}{g}
+ (2\pi \I)^2  dH_{\sigma_f}\wedge dH_{\sigma_g}.
\end{split}
\end{equation*}
The integral of $d\,\Log f \wedge d\,\Log g=d(\Log f \wedge d\,\Log g)$ over $M$ is zero because
$M$ is closed, and the integral of the last member,
$(2\pi \I)^2  dH_{\sigma_f}\wedge dH_{\sigma_g}$,
is an integer multiple of $(2\pi \I)^2$. Therefore, after integration
and taking the exponential we get
\begin{equation*}
\exp[\frac{1}{2\pi \I} \int_M \frac{df}{f}\wedge d\, \Log g
+\frac{1}{2\pi \I} \int_M d\,\Log f \wedge \frac{dg}{g} ]=1.
\end{equation*}
This proves the symmetry:
\begin{corollary}
Let $f$ and $g$ be two meromorphic functions on a closed Riemann
surface with disjoint divisors. Then
\begin{equation*}
\Res (f,g)=\Res (g,f).
\end{equation*}
\end{corollary}

\begin{remark}
This symmetry is also a consequence of Weil's reciprocity law \cite{Weil} (see Section~\ref{sec:local} for further details), and may alternatively be proved, in a more classical fashion, by
evaluating the integral in Cauchy's formula
$\int_{\partial M^\prime} \Log f \wedge d\,\Log g=0$ (cf. \cite[p.~242]{Griffith-Harris}). It is also  obtained by directly evaluating the last integral in (\ref{res2}).
\end{remark}

\begin{remark}
If the divisors of $f$ and $g$ are not disjoint but $f,g$ still form an admissible pair, then both $\Res (f,g)$ and $\Res (g,f)$ are either 0 or $\infty$, hence the symmetry remains valid although in a degenerate way. In this case, and more generally for nonadmissible pairs, Weil's reciprocity law in the form (\ref{Weil0}) (in Section~\ref{sec:local}) contains more information.
\end{remark}

By conjugating $g$ one gets the following formula for the modulus of
the resultant in terms of a Dirichlet integral.

\begin{theorem} \label{th:5}
Let $f$ and $g$ be two meromorphic functions on a compact Riemann
surface whose divisors have disjoint supports. Then
\begin{equation}\label{poten}
|\Res (f,g)|^2 = \exp[\frac{1}{2\pi \I}\int_M \frac{df}{f} \wedge \frac{d\bar{g}}{\bar{g}}].
\end{equation}
\end{theorem}


\begin{proof}
By Lemma~\ref{lem2} we have
$$
\frac{1}{2\pi \I} d\,\Log f \wedge d\,\Log \bar{g}
=\frac{1}{2\pi \I}\frac{df}{f} \wedge \frac{d\bar{g}}{\bar{g}}
+ \frac{df}{f} \wedge dH_{\sigma_g}
- dH_{\sigma_f} \wedge \frac{d\bar{g}}{\bar{g}}
-2\pi \I  dH_{\sigma_f} \wedge dH_{\sigma_g} .
$$
Integrating over $M$ and taking the exponential  yields, in view of (\ref{res2}), the required formula.
\qed\end{proof}

\section{Potential theoretic interpretations}
\label{sec5}

\subsection{The mutual energy and the resultant}
We recall some potential theoretic concepts (see, e.g., \cite{Saff} for more details). The potential of a signed measure (``charge distribution'') $\mu$ with compact support in $\Com{}$ is
\begin{equation*}
U^\mu(z)=-\int\log |z-\zeta|\;d\mu(\zeta).
\end{equation*}
The mutual energy between two such measures, $\mu$ and $\nu$, is (when defined)
\begin{equation*}\label{dis4}
\begin{split}
I(\mu,\nu)&=-\iint\log |z-\zeta|\;d\mu(z)d\nu(\zeta)=\int U^\mu\;d\nu=\int U^\nu\;d\mu,
\end{split}
\end{equation*}
and the energy of $\mu$ itself is $I(\mu)=I(\mu,\mu).$
In case $\int d\nu=\int d\mu=0$ the above mutual energy can after partial integration be written as a Dirichlet integral:
\begin{equation}\label{Dirichlet}
I(\mu,\nu)=
\frac{1}{2\pi}\int d U^\mu \wedge * d U^\nu,
\end{equation}
where $*$ is the Hodge star.

If $K\subset \Com{}$ is a compact set then either $I(\mu)=+\infty$ for all $\mu\geq 0$ with $\mathrm{supp}\,\mu\subset K$, $\int d\mu=1$, or
there is a unique such measure for which $I(\mu)$ has a finite minimum value. In the latter case $\mu$ is called the \textit{equilibrium distribution} for $K$ because its potential is constant on $K$ (except possibly for a small exceptional set), say
\begin{equation*}
U^\mu=\gamma \;\; (\mathrm{const}) \quad \text{on $K$}.
\end{equation*}
The logarithmic capacity of $K$ is defined as
\begin{equation*}
\mathrm{cap}\, (K)=e^{-\gamma}=e^{-I(\mu)}.
\end{equation*}
(If $I(\mu)=+\infty$ for all $\mu$ as above then $\mathrm{cap}\,(K)=0$).

Now let us think of signed measures as (special cases of) 2-form currents. Then, for example, (\ref{mud}) associates to each divisor $D$ in $\Com{}$ the charge distribution $\mu=\mu_D$. In particular, for any rational function $f$ of the form $f(z)=\prod_{i=1}^m\frac{z-a_i}{z-b_i}$ we have the charge distribution
\begin{equation*}
\mu=\mu_{(f)}=\sum_{i=1}^m \delta_{a_i}dx\wedge dy -\sum_{i=1}^m \delta_{b_i}dx\wedge dy ,
\end{equation*}
the potential of which is
$ 
U^\mu=-\log |f|.
$ 

One point we wish to make is that the resultant of two rational functions, $f$ and $g$, relates in the same way to the mutual energy.
In fact, with $\mu=\mu_{(f)}$ and $\nu=\mu_{(g)}$,
\begin{equation*}
\begin{split}
|\Res(f,g)|^2&=\exp [\scal{(f)}{\log g}+\scal{(f)}{\overline{\log g}}]
=e^{2\scal{(f)}{\log |g|}}=e^{-2\int U^\nu \,d\mu}=e^{-2I(\mu,\nu)},\\
\end{split}
\end{equation*}
hence
\begin{equation}
\label{nr}
I(\mu,\nu)=-\log |\Res(f,g)|. 
\end{equation}
The Dirichlet integral  (\ref{Dirichlet}) for $I(\mu,\nu)$ essentially gives the link between (\ref{nr}) and (\ref{poten}).

\subsection{Discriminant}
Recall that the (polynomial) discriminant $\mathrm{Dis}_{\mathrm{pol}}(f)$ is a polynomial in the coefficients of $f$ which vanishes whenever $f$ has a multiple root.
In case of a \textit{monic} polynomial $f(z)=\prod_{i=1}^m (z-a_i)$ we have
\begin{equation*}
\mathrm{Dis}_{\mathrm{pol}}(f)=(-1)^{\frac{m(m-1)}{2}}\Res_{\mathrm{pol}}(f,f')=\prod_{i< j} (a_i-a_j)^2.
\end{equation*}
Thus the discriminant is the square of the Van der Monde determinant.

The discriminant can be related to a renormalized self-energy of the measure $\mu=\mu_{(f)}$. The self-energy itself is actually infinite because point charges always have infinite energy. Formally:
\begin{equation*}
\begin{split}
I(\mu)&=\int U^\mu \,d\mu=\scal{(f)}{-\log |f|}=-\log \prod_{i,j=1}^m |a_i-a_j| \quad (=+\infty).
\end{split}
\end{equation*}
The renormalized energy $\widehat I(\mu)$ is obtained by simply subtracting off the infinities $I(\delta_{a_i})$, i.e., the diagonal terms above:
\begin{equation*}
\begin{split}
\widehat{I}(\mu)&=-\log \prod_{i\ne j} |a_i-a_j| =-\log \prod_{i< j} |a_i-a_j|^2=-\log |\mathrm{Dis}_{\mathrm{pol}}(f)|.
\end{split}
\end{equation*}

Thus,
$ 
|\mathrm{Dis}_{\mathrm{pol}}(f)|=e^{-\widehat{I}(\mu)}.
$ 
Here $\int d\mu=\deg f=m$, and after normalization (there are $m(m-1)$ factors in $\mathrm{Dis}_{\mathrm{pol}}(f)$)
it is known that the transfinite diameter
\begin{equation*}
d_{\infty}(K)=\lim_{m\to \infty}\; \max_{\deg f=m} |\mathrm{Dis}_{\mathrm{pol}}(f)|^{\frac{1}{m(m-1)}},
\end{equation*}
equals the capacity:
$ 
\label{capacity}
d_{\infty}(K)=\mathrm{cap}\,(K).
$ 

Notice also that the discriminant may be regarded as a renormalized \textit{self-resultant} $\Res_{\mathrm{pol}}(f,f)$:
\begin{equation}\label{dis_renorm}
\Res_{\mathrm{pol}}(f,f)=\prod_{i,j} (a_i-a_j)\;\stackrel{\text{renorm}}{\Longrightarrow}
\;\mathrm{Dis}_{\mathrm{pol}}(f)=\prod_{i\ne j} (a_i-a_j).
\end{equation}

We can use the same renormalization method to arrive at a definition of discriminant in the rational case. Let $f$ be a rational function
\begin{equation*}
f(z)=\frac{f_1(z)}{f_2(z)}\equiv  \frac{\prod_{i=1}^m(z-a_i)}{\prod_{i=1}^m(z-b_i)}.
\end{equation*}
Then applying the scheme  in (\ref{dis_renorm}) gives
\begin{equation}\label{dis_renorm1}
\begin{split}
&\Res(f,f)=\prod_{i,j} \frac{(a_i-a_j)(b_i-b_j)}{(a_i-b_j)(b_i-a_j)}\quad \stackrel{\text{renorm}}{\Longrightarrow}\\
\stackrel{\text{renorm}}{\Longrightarrow}\quad
&\mathrm{Dis}(f):=\frac{\prod\limits_{i\ne j}(a_i-a_j)\prod\limits_{i\ne j}(b_i-b_j)}{\prod\limits_{i,j}(a_i-b_j)\prod\limits_{i,j}(b_i-a_j)}
=\frac{\Res_{\mathrm{pol}}(f_1,f_1')\Res_{\mathrm{pol}}(f_2,f_2')}{\Res_{\mathrm{pol}}(f_1,f_2)\Res_{\mathrm{pol}}(f_2,f_1)}.
\end{split}
\end{equation}

The corresponding renormalized energy of $\mu=\mu_{(f)}$ is
\begin{equation*}
\begin{split}
\widehat{I}(\mu)&=-\log \left|\frac{\prod_{i\ne j}(a_i-a_j)\prod_{i\ne j}(b_i-b_j)}{\prod_{i,j}(a_i-b_j)\prod_{i,j}(b_i-a_j)}\right|=-\log |\mathrm{Dis}(f)|
\end{split}
\end{equation*}
which yields
\begin{equation*}
|\mathrm{Dis}(f)|=e^{-\widehat{I}(\mu)}.
\end{equation*}

We note that the definition (\ref{dis_renorm1}) of $\mathrm{Dis}(f)$ is consistent with the so-called characteristic property of the polynomial discriminant \cite[p.~405]{Gelfand-Kapranov-Zelevinskij}. Namely, one can easily verify that the meromorphic resultant of two rational functions can be obtained as the polarization of the discriminant in (\ref{dis_renorm1}), that is
\begin{equation*}
\Res(f,g)^2=\frac{\mathrm{Dis}(fg)}{\mathrm{Dis}(f)\mathrm{Dis}(g)}.
\end{equation*}

\subsection{Riemann surface case}
Much of the above can be repeated for an arbitrary compact Riemann surface $M$. For any signed measure $\mu$ on $M$ with $\int_Md\mu=0$ there is potential $U^\mu$, uniquely defined up to an additive constant, such that
\begin{equation*}
-d * d U^\mu=2\pi \mu.
\end{equation*}
Here $\mu$ is considered as a 2-form current ($\mu$ may actually be an arbitrary 2-form current with $\scal{\mu}{1}=0$, and then $U^\mu$ will be a 0-current; the existence and uniqueness of $U^\mu$ follows from ordinary Hodge theory, see e.g. \cite[p.~92]{Griffith-Harris}).

The mutual energy between two measures as above can still be defined as
\begin{equation*}
\begin{split}
I(\mu,\nu)&=\int U^\mu\;d\nu=\int U^\nu\;d\mu
\end{split}
\end{equation*}
and (\ref{Dirichlet}) remains true. Similarly, (\ref{nr}) remains valid for $\mu=\mu_{(f)}$, $\nu=\mu_{(g)}$. Thus
\begin{equation*}
|\Res(f,g)|=e^{-I(\mu,\nu)}.
\end{equation*}

It is interesting to notice that this gives a way of defining the modulus of the resultant of any two divisors of degree zero: if $\deg D_1=\deg D_2=0$ with $\supp D_1 \cap \supp D_2=\emptyset$ then one naturally sets
\begin{equation*}
|\Res(D_1,D_2)|=e^{-I(\mu_{D_1},\mu_{D_2})}.
\end{equation*}
It is not clear whether there is any natural definition of $\Res(D_1,D_2)$ itself, except in genus zero where we have (\ref{crossr}). Directly from the definition (\ref{main1new}) we can however define $\Res(D,g)=g(D)$ for $D$ a divisor of degree zero and $g$ a meromorphic function.

\section{The resultant as a function of the quotient}
\label{sec:sze}

\subsection{Resultant identities}\label{seq:residen}

In previous sections we have considered the resultant as a function of two meromorphic functions, $f$ and $g$, say.
Sometimes, however, it is possible and convenient to think of the resultant as a function of just \textit{one} function, namely the quotient $h=\frac{f}{g}$. In general, part of the information about $f$ and $g$ is lost in $h$, hence some additional information has to be provided.

For instance, if $f$ and $g$ are two \textit{monic} polynomials, then formula (\ref{hankel}) in its simplest form, when $N=n$, reads
\begin{equation*}
\Res_{\mathrm{pol}}(f,g)=\det t_{m,n}(h).
\end{equation*}

Another example is if the divisors of $f$ and $g$ are confined to lie in prescribed disjoint sets: given any set $U\subset M$ then among pairs $f,g$ with $\supp (f) \subset U$, $\supp (g)\subset M\setminus U$, the resultant $\Res(f,g)$ only depends on $\frac{f}{g}$.
Integral representations for $\Res(f,g)$ in terms of only $f/g$ and $U$ will in such cases be elaborated in Section~\ref{subsec:63} (Theorem~\ref{the:coh}).

In the remaining part of this section we shall pursue a further point of view. Suppose that the divisors of $f$ and $g$ are not necessarily disjoint but that $f$ and $g$ still form an admissible pair. In general we have, with $h=f/g$,
\begin{equation*}\label{problord}
\ord h\leq \ord f+\ord g,
\end{equation*}
and it is easy to see that $\Res(f,g)=0$ if and only if this inequality is strict (because strict inequality means that at least one common zero or one common pole of $f$, $g$ cancels out in the quotient $f/g$).

Now start with $h$ and consider admissible pairs $f,g$ with $h=f/g$ and such that
\begin{equation}\label{problord1}
\ord h= \ord f+\ord g.
\end{equation}
In general there are many such pairs $f,g$ and by the above $\Res(f,g)\ne 0$ for all of them. The question we want to consider is whether there are any restrictions on which values $\Res(f,g)$ can take. At least in the rational case there turns out be such restrictions and this is what we call \textit{resultant identities}.

Let $d\geq 1$ and
\begin{equation}\label{hh}
h(z)= \prod_{i=1}^d\frac{z-a_i}{z-b_i}.
\end{equation}
Let $C_d^m$ denote the set of all increasing length-$m$ sequences $(i_1,\ldots,i_m)$, $1\le i_1<\ldots<i_m\le d$. For two given elements $I,J\in C_d^m$ define
\begin{equation*}
h_{IJ}(z)=\frac{\prod_{i\in I}(z-a_i)}{\prod_{j\in J}(z-b_j)},
\end{equation*}
Then all the solutions $f$, $g$ of (\ref{problord1}), up to a constant factor (which by (\ref{homm}) is inessential for the resultant), are parameterized  by
\begin{equation*}\label{anal}
f(z)=h_{IJ}(z),
\qquad g(z)=\frac{h_{IJ}(z)}{h(z)}=\frac{1}{h_{I'J'}(z)},
\end{equation*}
where the prime denotes complement, e.g., $I'=\{1,\ldots,d\}\setminus I$.

The main observation of this section is that the resultants $\Res(f,g)$ satisfy a system of linear identities. An extended version of the material below with applications to rational and trigonometric identities will appear in \cite{GT071}.

\begin{proposition}\label{the:o}
Let $0\leq m\leq d$ and $J\in C_d^m$. Then
\begin{equation}\label{res:sup}
\sum_{I\in C_d^m}\Res(h_{IJ}, 1/h_{I'J'})=\sum_{I\in C_d^m}\Res(h_{JI}, 1/h_{J'I'})=1.
\end{equation}
\end{proposition}

\begin{proof}
We briefly describe the idea of the proof.
Denote by $\mathbf{A}$ and $\mathbf{B}$ the two Van der Monde matrices with entries $(a_i^{j-1})$ and $(b_i^{j-1})$, $1\leq i,j\leq d$, respectively. Let $I=\{{i_1},\ldots,{i_m}\}$ and $J=\{{j_1},\ldots,{j_m}\}$. Then one can readily show that
\begin{equation}\label{res:sup3}
\Res(h_{IJ}, 1/h_{I'J'})=(-1)^n\det \Lambda_{IJ} \det (\Lambda^{-1})_{IJ},
\end{equation}
where $n={\sum_{s=1}^m (i_s+j_s)}$.
Here $\Lambda=\mathbf{A}\mathbf{B}^{-1}$ and $\Lambda_{IJ}$ (resp. $(\Lambda^{-1})_{IJ}$) denotes the minor of $\Lambda$ (resp. $\Lambda^{-1}$)
formed by intersection of the rows $i\in I$ and the columns $j\in J$. Hence the required identities follow from (\ref{res:sup3}) and the Laplace expansion theorem for determinants.
\qed\end{proof}

In the simplest case, $d=2$, $m=1$, (\ref{res:sup}) amounts to the characteristic property of the cross-ratio:
\begin{equation*}
(a,b,c,d)+(a,c,b,d)=1.
\end{equation*}

The resultants in (\ref{res:sup}) appear also in the so-called Day's formula \cite{day} for the determinants of truncated Toeplitz operators. Let $h$ be a function  given by (\ref{hh}) such that $|b_i|\ne1$ for all $i$, and let
$ J=\{j:\; |b_j|>1\}$.

Introduce the Toeplitz matrix of order $N$
\begin{equation}\label{days}
t_N(h)\equiv \left(
              \begin{array}{cccc}
                h_{0} & h_{-1} &  \ldots & h_{1-N}\\
                h_{1} & h_{0} & \ldots & h_{2-N} \\
                \ldots & \ldots & \ldots & \ldots \\
                h_{N-1} & h_{N-2} &\ldots & h_0
              \end{array}
            \right)
\end{equation}
where $h_k=\frac{1}{2\pi }\int_{0}^{2\pi}e^{-\I k\theta}h(e^{\I \theta})d\theta$ are the Fourier coefficients of $h$ on the unit circle.
Then, in our notation, Day's formula reads
\begin{equation}\label{days1}
\det t_N(h)=\sum_{I\in C_d^m}\Res(h_{IJ}, 1/h_{I'J'})\cdot h^N_{I'J'}(0),
\end{equation}
where $m$ denotes the cardinality of $J$ and $N\geq 1$. Notice that formal substitution of $N=0$ with $t_0(h)=1$ into (\ref{days1}) gives exactly  the statement of Proposition~\ref{the:o}.

\begin{remark}
Taking  double sums in (\ref{res:sup}) (over all $I,J\in C_d^m$) we get quantities which occur also when computing subresultants (see, e.g., \cite{LascPr}). Recall that the (scalar) subresultant of degree $k$ is the determinant of the matrix obtained from the Sylvester matrix (\ref{sfor}) by deleting the last $2k$ rows and the last $k$ columns with coefficients of $f$, and the last $k$ columns with coefficients of $g$. In a different context, the subresultants are determinants of certain submatrices of the Sylvester matrix (\ref{sfor}) which occur as successive remainders in finding the greatest common divisor of two polynomials by the Euclid algorithm \cite{Syl40}.
\end{remark}

The identities (\ref{res:sup}) have beautiful trigonometric interpretations. Take
\begin{equation*}
f(z)=\prod_{k=1}^m\frac{z-e^{2\I a_k}}{z-e^{2\I b_k}},\qquad
g(z)=\prod_{l=1}^n\frac{z-e^{2\I c_l}}{z-e^{2\I d_l}}.
\end{equation*}
Then one easily finds that
\begin{equation*}\label{trig10}
\Res(f,g)=\prod_{k=1}^m\prod_{l=1}^n
\frac{\sin(a_k-c_l)}{\sin(a_k-d_l)}\frac{\sin(b_k-d_l)}{\sin(b_k-c_l)},
\end{equation*}
hence a direct application of (\ref{res:sup}) gives the following.

\begin{corollary}\label{cor:dress}
Let $d\geq 2$ and $J\in C_d^m$. Then
\begin{equation}\label{res:sup0}
\sum_{I}\frac{\prod_{i,j'}\sin(a_i-b_{j'})\prod_{i',j}\sin(b_{j}-a_{i'})}{\prod_{i,i'}\sin(a_i-a_{i'})\prod_{j,j'}\sin(b_{j}-b_j')}=1,
\end{equation}
where the sum is taken over all subsets $I\in C_d^m$ and the product over $i\in I$, $i'\in I'$, $j\in J$, $j'\in J'$.
\end{corollary}

For example, specializing by taking $b_j=\frac{\pi}{2}+a_i$  in (\ref{res:sup0}) one gets identities in the spirit of those given recently in \cite{Calog}, \cite{Calog1}.

There are also analogues of Proposition~\ref{the:o} for the complex torus $M=\Com{}/L_\tau$. For these  one has to take into account the Abel condition (\ref{jacobcond}). Although we have not been able to find complete analogues of the rational resultant identities, one particular case is worth mentioning here. Notice that  the minimal possible value of $d$ in order for a meromorphic function $h(z)=\prod_{i=1}^d\frac{ \theta (z-u_i)}{\theta (z-v_i)}$ to split into two \textit{non-constant} meromorphic functions, i.e. $h=f/g$, is $d=4$. One can readily show that any such function may be written as
\begin{equation*}
h(z)=\frac{\phi(z-z_0,a_1)\phi(z-z_0,a_2)}{\phi(z-z_0,b_1)\phi(z-z_0,b_2)},
\end{equation*}
where $\phi(\zeta,a)=\theta (\zeta-a)\theta (\zeta+a)$. We additionally assume that $a_1\pm a_2\not\in L$ and $b_1\pm b_2\not\in L$.
Then all non-constant solutions of (\ref{problord1}) are given by
\begin{equation*}
f(z)=\frac{\phi(z,{a_i})}{\phi(z,b_j)},\quad g(z)=\frac{\phi(z,{b_{j'}})}{\phi(z,{a_{i'}})},  \qquad i,j=1,2,
\end{equation*}
where $\{k,k'\}=\{1,2\}$. Hence
\begin{equation*}
\rho_{ij}:=  \Res(f,g)=\left[\frac{\theta(a_i-b_{j'})\theta(a_i+b_{j'})\theta(a_{i'}-b_j)\theta(a_{i'}+b_j)}
{\theta(a_i-a_{i'})\theta(a_i+a_{i'})\theta(b_j-b_{j'})\theta(b_j+b_{j'})}\right]^2,
\end{equation*}
and there only two different values of $\rho_{ij}$:
\begin{equation*}
\xi_1:=\rho_{11}=\rho_{22},\qquad \xi_2:=\rho_{12}=\rho_{21}.
\end{equation*}
Using the famous addition theorem of Weierstra{\ss}
\begin{equation*}\label{jacobI}
\begin{split}
0=&\theta(a-c)\theta(a+c)\theta(b-d)\theta(b+d)-\theta(a-b)\theta(a+b)\theta(c-d)\theta(c+d)\\
-&\theta(a-d)\theta(a+d)\theta(b-c)\theta(b+c),\\
\end{split}
\end{equation*}
one  finds that (with appropriate choices of signs)
\begin{equation}\label{likei}
\pm\sqrt{\xi_1}\pm\sqrt{\xi_2}=1,
\end{equation}
or more adequately:
$
(1-\xi_1)^2+(1-\xi_2)^2=2\xi_1\xi_2.
$

The identity (\ref{likei}) may be generalized to functions of the kind
\begin{equation*}
h(z)=\prod_{k=1}^d\frac{\phi(z-z_0,a_k)}{\phi(z-z_0,b_k)}.
\end{equation*}
However the problem of description of the range of $\Res(f,g)$ in (\ref{problord1}) for general meromorphic functions $h$ on $\Com{}/L_\tau$ remains  open.

\subsection{Integral representation of $\Res_{U}$}
\label{subsec:63}

Let us now turn to the situation of having a preassigned set $U\subset
M$ and consider resultants $\Res(f,g)$ for meromorphic functions
$f$ and $g$ with $\supp (f)\subset U$, $\supp (g)\subset M\setminus U$.
It is easy to see that for such pairs $\Res(f,g)$ only depends
on the quotient $h=f/g$. Indeed, this is obvious from the fact (see
(\ref{homm})) that the resultant only depends on the divisors:
under the above assumptions the divisors of $f$ and $g$ are clearly
determined by $h$ and $U$.

To make the above in a slightly more formal we may define
$\Res(D_1,D_2)$ for any two principal divisors $D_1$, $D_2$ having, e.g.,
disjoint supports. For any divisor $D$, let $D_U$ denote its restriction
to the set $U$ and extended by zero outside $U$
(thus with $D=\sum_{a\in M} D(a)a$, $D_U=\sum_{a\in U} D(a)a$).
Then in the situation at hand we can write
\begin{equation*}
\Res(f,g)=\Res((f),(g))=\Res((h)_U,(h)_U-(h)),
\end{equation*}
which only depends on $h$ and $U$. This motivates the following definition.

\begin{definition}
For any set $U\subset M$ and any meromorphic function $h$ on $M$ such
that $(h)_U$ is a principal divisor we define
\begin{equation*}
\Res_U(h)=\Res((h)_U,(h)_U-(h)).
\end{equation*}
\end{definition}

It is easy to check that
\begin{equation*}
\Res_{U}(h)=\Res_{M\setminus U}(h).
\end{equation*}
We shall consider the symmetric situation that
\begin{equation*}
M=U\cup \Gamma \cup V,
\end{equation*}
where $U$, $V$ are disjoint nonempty open sets and $\Gamma=\partial
U=\partial V$. We provide $\Gamma$ with the orientation of $\partial U$. By the above, with $f$ and $g$ meromorphic on $M$, $\supp (f)\subset U$, $\supp (g)\subset V$ and $h=f/g$ we have
\begin{equation*}
\Res_U(h)=\Res_V(h)=\Res(f,g).
\end{equation*}
Note that the function $h$ is holomorphic and nonzero in a neighbourhood of $\Gamma$, $h\in \mathcal{O}^*(\Gamma)$, and that it is uniquely defined by its values on $\Gamma$. Our aim is to find an integral representation for $\Res_U(h)$ in terms only of the values of $h$ on $\Gamma$.

The problem of decomposing a given $h\in \mathcal{O}^*(\Gamma)$ into functions $f\in \mathcal{O}^*(\overline{V})$, $g\in \mathcal{O}^*(\overline{U})$ with $h=f/g$ is a special case of the second Cousin problem. By taking logarithms we shall reduce it, under symplifying assumptions, to the corresponding additive problem, which is the first Cousin problem. For the latter we have the following simple criterion for solvability.

\begin{lemma}
\label{lemA}
Let $M=U\cup \Gamma\cup V$ be as above. Necessary and sufficient condition for a function $H\in \mathcal{O}(\Gamma)$ to be  decomposable as
\begin{equation*}
H=H_+-H_- \qquad \text{on $\Gamma$}
\end{equation*}
with $H_+\in \mathcal{O}(\overline{U})$, $H_-\in \mathcal{O}(\overline{V})$ is that
\begin{equation*}\label{Homega}
\int_\Gamma H\wedge \omega=0 \qquad \text{for all $\omega\in \mathcal{O}^{1,0}(M)$}.
\end{equation*}
When the decomposition exists the functions $H_\pm$ are unique up to addition of a common constant (more adequately: a function in $\mathcal{O}(M)$).
\end{lemma}

The lemma is well-known and can be deduced for example from the Serre duality theorem. We shall just remark that ``explicit'' representations of $H_\pm$ can be given in terms of a suitable Cauchy kernel:
\begin{equation*}
H_{\pm}(z)= \frac{1}{2\pi\I}\int_\Gamma H(\zeta)\Phi(z,\zeta;z_0,
\zeta_0)\,d\zeta
\end{equation*}
the plus sign for $z\in U$, minus for $z\in V$. The kernel $\Phi(z,\zeta;z_0,
\zeta_0)$ is, in the variable $z$, a meromorphic function with a simple pole at $z=\zeta$
and a pole of higher order (depending on the genus) at
$z=\zeta_0$. In the variable $\zeta$ it is a meromorphic one-form
with simple poles of residues plus and minus one at $\zeta=z$ and
$\zeta=z_0$ respectively; $z_0$ and $\zeta_0$ are fixed but
arbitrary points, $z_0\ne \zeta_0$. In the case of the Riemann
sphere, $\Phi(z,\zeta;z_0, \zeta_0)\,d\zeta$ is the ordinary Cauchy
kernel
\begin{equation}\label{Cauchykernel}
\Phi(z,\zeta;z_0, \zeta_0)\,d\zeta =\frac{d\zeta}{\zeta-z}
-\frac{d\zeta}{\zeta-z_0},
\end{equation}
hence does not involve $\zeta_0$. In the the case of higher genus
the point $\zeta_0$ is really needed. We refer to \cite{Rodin} for
the construction of the Cauchy kernel in general.

\begin{theorem}\label{the:coh}
Let $M=U\cup \Gamma \cup V$ with $U$ connected and simply connected, and let $h$ be meromorphic on $M$ without poles and zeros on $\Gamma$. Assume in addition that
\begin{equation}\label{dhh}
\frac{1}{2\pi \I}\int_\Gamma \frac{dh}{h}=0
\end{equation}
and that
\begin{equation}\label{Loghomega}
\int_\Gamma\Log h \wedge \omega =0 \quad \text{for all $\omega\in
\mathcal{O}^{1,0}(M)$}
\end{equation}
(the previous condition guarantees that a single-valued branch of $\log h$ exists on $\Gamma$). Then $(h)_U$ is a principal divisor and
\begin{equation*}
\Res_U(h)=\exp \;[\,\frac{1}{2\pi \I}\int_\Gamma d\,(\Log
h)_-\wedge(\Log h)_+].
\end{equation*}
\end{theorem}

\begin{remark}
Ideally (\ref{Loghomega}) should be replaced be the weaker condition that there exists a closed 1-chain $\gamma$ on $M$ such that
\begin{equation}\label{Abel}
\int_\Gamma \Log h\wedge \omega=2\pi \I \int_\gamma \omega \qquad \text{for all $\omega\in\mathcal{O}^{1,0}(M)$}.
\end{equation}
In fact, this turns out to be exactly, by Abel's theorem, the necessary and sufficient condition for $(h)_U$ to be a principal divisor. However, (\ref{Abel}) would lead to a more complicated formula for $\Res_U(h)$. Note that (\ref{Loghomega}) is vacuously satisfied in the case $M=\mathbb{P}$, which will be our main application. Condition (\ref{dhh}) says that the divisor $(h)_U$ has degree zero.
\end{remark}

\begin{proof}
We first prove that $(h)_U$ is a principal divisor. Using the notation of Lemma~\ref{lem2} we make $\Log h$ into a single-valued function on all of $M$ by making cuts along a 1-chain $\sigma_h$ such that $\partial \sigma_h=(h)$. Since $\Log h$ is already single-valued on $\Gamma$, $\sigma_h$ can be chosen not to intersect $\Gamma$. Thus $\sigma_h$ consists of two disjoint parts, $\sigma_h\cap U$ and $\sigma_h \cap V$. The terms of $\sigma_h$ containing the curves $\alpha_1,\ldots,\beta_p$ will appear in $\sigma_h\cap V$ because $U$ is simply connected.

Now, for all $\omega\in\mathcal{O}^{1,0}(M)$ we have by (\ref{Loghomega}) and Lemma~\ref{lem2}
\begin{equation*}
\begin{split}
0&=\frac{1}{2\pi \I} \int_\Gamma \Log h\wedge \omega=\frac{1}{2\pi \I} \int_U d\Log h\wedge \omega
=\frac{1}{2\pi \I} \int_U \left(\frac{dh}{h}-2\pi \I d H_{\sigma_h}\right)\wedge \omega\\
&=\frac{1}{2\pi \I} \int_U \frac{dh}{h}\wedge \omega -\int_U d H_{\sigma_h}\wedge \omega
=-\int_M d H_{\sigma_h\cap U}\wedge \omega=-\int_{\sigma_h\cap U}\omega.
\end{split}
\end{equation*}

By Abel's theorem this implies that $\partial(\sigma_h\cap U)=(h)_U$  is a principal divisor (condition (\ref{Abel}), in place of (\ref{Loghomega}), would have been enough for this conclusion).

The divisor $(h)_U$ being principal means that $(h)_U=(f)$ for some $f$ meromorphic on $M$. Setting $g=f/h$ we have $\supp (f)\subset U$, $\supp (g)\subset V$ and $h=f/g$. It follows that $\Res_U(h)=\Res(f,g)$, hence to prove the theorem it is by Theorem~\ref{integralformula} enough to prove that
$$
\int_\Gamma d\,(\Log
h)_-\wedge(\Log h)_+=\int_M \frac{df}{f}\wedge d \, {\Log \,} g.
$$

To that end we shall compare two decompositions of $d\Log h=\frac{dh}{h}$ on $\Gamma$: from Lemma~\ref{lemA} we get
$$
d\Log h=d(\Log h)_+-d(\Log h)_- \qquad \text{on $\Gamma$}
$$
with $(\Log h)_+\in \mathcal{O}(\overline{U})$, $(\Log h)_-\in \mathcal{O}(\overline{V})$, while $h=f/g$ gives
$$
\frac{dh}{h}=\frac{df}{f}-\frac{dg}{g} \qquad \text{on $\Gamma$},
$$
where $df/f\in \mathcal{O}^{1,0}(\overline{V})$, $dg/g\in \mathcal{O}^{1,0}(\overline{U})$.

It follows that
$$
\frac{df}{f}+d(\Log h)_-=\frac{dg}{g}+d(\Log h)_+  \qquad \text{on $\Gamma$}
$$
and that the left and right members combine into a global 1-form $\omega_0\in \mathcal{O}^{1,0}(M)$. Thus
\begin{equation*}
\begin{split}
d(\Log h)_-&=\omega_0-\frac{df}{f} \;\; \text{in $V$},\qquad
d(\Log h)_+=\omega_0-\frac{dg}{g} \;\; \text{in $U$}.\\
\end{split}
\end{equation*}

In the simply connected domain $U$ we may write $\omega_0=d\varphi$ for some $\varphi\in \mathcal{O}(\overline{U})$ and also
$
\frac{dg}{g}=d\,\Log g
$
($dH_{\sigma_g}=0$ in $U$ because $\sigma_g$ can be chosen to be $\sigma_h\cap V$; similarly $\sigma_f$ can be chosen to be $\sigma_h\cap U$).
It follows after integration and adjusting $\varphi$ by a constant that
$$
(\Log h)_+=\varphi-\Log g \qquad \text{in $U$}.
$$

Now we finally obtain
\begin{equation*}
\begin{split}
&\int_\Gamma d\,(\Log
h)_-\wedge(\Log h)_+=\int_\Gamma (\omega_0-\frac{df}{f})\wedge(\varphi-\Log g )
=-\int_\Gamma \frac{df}{f}\wedge(\varphi-\Log g )\\
&=\int_{V} \frac{df}{f}\wedge d\Log g -\int_{\Gamma}(d\Log h+d\Log g)\wedge \varphi
=\int_{M} \frac{df}{f}\wedge d\Log g ,
\end{split}
\end{equation*}
as desired.
\qed\end{proof}

\begin{remark}
Under the assumptions of the theorem, the solution of the second Cousin problem of finding $f,g$ such that $h=f/g$ on $\Gamma$ is given by
\begin{equation*}
\begin{split}
f&=\exp\biggl[\int \frac{df}{f}\biggr] =\exp\biggl[\int (\omega-d(\Log h)_-)\biggr]\qquad \text{in $V$},\\
g&=\exp\biggl[\int \frac{dg}{g}\biggr] =\exp\biggl[\int (\omega-d(\Log h)_+)\biggr]\qquad \text{in $U$}\\
\end{split}
\end{equation*}
(indefinite integrals), where $\omega\in \mathcal{O}^{1,0}(M)$ is to be chosen such that $\int (\omega-d(\Log h)_-)$ is single-valued in $V$ modulo multiples of $2\pi \I$.
\end{remark}

\subsection{Cohomological interpretations of the quotient}\label{linebundle}
Let us give some interpretations of the above material in terms of \v{C}ech cohomology.
Given $h\in \mathcal{O}^*(\Gamma)$, let $U_1$, $V_1$ be open
neighbourhoods of $\overline{U}$ and $\overline{V}$, respectively,
such that $h\in \mathcal{O}^*(U_1\cap V_1)$. Then $\{U_1,V_1\}$ is
an open covering of $M$, and relative to this $h$ represents an
element $[h]$ in $H^1(M,\mathcal{O}^*)$. It is well-known
\cite{Gunning72}, \cite{Forster} that $[h]=0$ as an element in
$H^1(M,\mathcal{O}^*)$ if and only if $h$ is a coboundary already
with respect to $\{U_1,V_1\}$, i.e., if and only if there exist
$f\in\mathcal{O}^*(V_1)$ and $g\in\mathcal{O}^*(U_1)$ such that $h=f/g$
in $U_1\cap V_1$. If $h$ is meromorphic in $M$, then so are $f$ and
$g$.

Similarly, a function $H\in \mathcal{O}(\Gamma)$ represents an element $[H]$ in $H^1(M,\mathcal{O})$, and $[H]=0$ if and only if there exist $F\in \mathcal{O}(U_1)$, $G\in \mathcal{O}(V_1)$ (for some $U_1\supset \overline{U}$, $V_1\supset \overline{V}$) such that $H=F-G$ on $\Gamma$.

The spaces $H^1(M,\mathcal{O})$ and $H^1(M,\mathcal{O}^*)$ are related via the long exact sequence of cohomology groups which comes from the exponential map on the sheaf level: with $e(f)=\exp [2\pi \I f]$ we have
$$
0\;  \rightarrow \; \mathbb{Z} \;  \rightarrow \; \mathcal{O} \;\stackrel{e}{\rightarrow}\; \mathcal{O}^* \;\rightarrow\; 1,
$$
hence
\begin{equation*}
\begin{split}
0  &\rightarrow  H^0(M,\mathbb{Z}) \rightarrow H^0(M,\mathcal{O}) \rightarrow
H^0(M,\mathcal{O}^*) \rightarrow H^1(M,\mathbb{Z}) \rightarrow \\
&  \rightarrow H^1(M,\mathcal{O}) \stackrel{e}{\rightarrow}
H^1(M,\mathcal{O}^*)  \rightarrow H^2(M,\mathbb{Z})   \rightarrow 0.\\
\end{split}
\end{equation*}

From this  we extract the exact sequence
\begin{equation}\label{exact}
0\;  \rightarrow \; {H^1(M,\mathcal{O})}/{H^1(M,\mathbb{Z})} \;  \stackrel{e}{\rightarrow} \; H^1(M,\mathcal{O}^*) \;\stackrel{c}{\rightarrow}\;
H^2(M,\mathbb{Z}) \;  \rightarrow \; 0.\\
\end{equation}
Here $c$ is the map which associates to $[h]\in H^1(M,\mathcal{O}^*)$  its characteristic class, or Chern class, and it is readily verified that it is given by
\begin{equation*}
\begin{split}
c([h])&=\wind_\Gamma h=\frac{1}{2\pi \I}\int_\Gamma \frac{dh}{h}=\deg (h)_U.
\end{split}
\end{equation*}

If $c([h])=0$ then $[h]$ is in the range of $e$. If $\Gamma$ is connected then $\log h$ is single-valued on $\Gamma$ and the preimage of $[h]$ can be represented by $H=\frac{1}{2\pi \I}\,\Log h$. However, if $\Gamma$ is not connected then the preimage of $[h]$ cannot always be represented by a function $H\in \mathcal{O}(\Gamma)$, one needs a finer covering of $M$ than $\{U_1,V_1\}$ to represent it. This is a drawback of the method using the decomposition $M=U\cup \Gamma\cup V$ in combination with the $\exp$--$\log$ map and explains some of our extra assumptions in Theorem~\ref{the:coh}.

Assume nevertheless that the preimage of $[h]\in H^1(M,\mathcal{O}^*)$ (with $c([h])=0$) can be represented by $H=\frac{1}{2\pi \I}\,\Log h\in \mathcal{O}(\Gamma)$. Then of course $[h]=0$ if $[H]=0$ as an element in $H^1(M,\mathcal{O})$, i.e., if $\int_\Gamma H\wedge \omega=0$ for all $\omega\in \mathcal{O}^{1,0}(M)$. However, what exactly is needed for $[h]=0$ is by (\ref{exact}) only that $[H]\in H^1(M,\mathbb{Z})$, and this what is expressed in (\ref{Abel}).

Since, for $H\in \mathcal{O}(\Gamma)$, $[H]=0$ as an element in $H^1(M,\mathcal{O})$ if and only if $\int_\Gamma H\wedge \omega=0$ for all $\omega\in \mathcal{O}^{1,0}(M)$, the pairing
\begin{equation*}
(\omega,H) \; \mapsto \; \int_\Gamma H\wedge \omega
\end{equation*}
descends to a bilinear map
\begin{equation*}
H^0(M,\mathcal{O}^{1,0})\times H^1(M,\mathcal{O}) \rightarrow \Com{}.
\end{equation*}
This map is in fact the Serre duality pairing (\cite{Serre 55}, \cite{Gunning72}) with respect to the covering $\{U_1,V_1\}$. Versions of the Serre duality with respect to more general coverings will be discussed in the next section.


\subsection{Resultant via Serre duality}\label{subsec:serre}

We now return to the general integral formula
in Theorem~\ref{integralformula}, and interpret the exponent $\frac{1}{2\pi \I}\int_M \frac{df}{f}\wedge d \, {\Log \,} g$ directly in terms of the Serre
duality pairing, which in general also involves a line bundle or a divisor. With a divisor $D$, the pairing looks
$$
\scal{\;}{\;}_{\mathrm{Serre}}:\;H^0(M,\mathcal{O}^{1,0}_D)\times H^1(M,\mathcal{O}_{-D})\to \C,
$$
between meromorphic $(1,0)$-forms with
divisor $\geq -D$ and (equivalence classes of) cocycles of
meromorphic functions with divisor $\geq D$.

In our case, given two meromorphic functions $f$ and $g$, we choose
$D\geq 0$ to be the divisor of poles of $\frac{df}{f}$ (or any
larger divisor), so that $\frac{df}{f}\in \Gamma(M,
\mathcal{O}_{D}^{1,0})$. As for the other factor, $\log g$ defines
an element, which we denote by $[\delta \log g]$, of $H^1 (M,
\mathcal{O}_{-D})$ as follows. First, with $\gamma_g$ as in the beginning of Section~\ref{sec:begin}, choose an open cover $\{U_i\}$
of $M$ consisting of simply connected domains $U_i$ satisfying
\begin{equation*}
(\supp D \cup \supp \gamma_g) \cap U_i \cap U_j = \emptyset \quad
\text{whenever\,\,} i\ne j
\end{equation*}
(in particular $\supp \gamma_g \cap \partial U_i =\emptyset$ for all
$i$). Second, choose for each $i$ a branch, $(\log g)_i$, of $\log
g$ in $U_i\setminus \gamma_g$.  Finally, define a cocycle $\{(\delta\log g)_{ij}\}$, to
represent $[\delta \log g]\in H^1(M, \mathcal{O}_{-D})$, by
\begin{equation*}
(\delta\log g)_{ij} = (\log g)_i - (\log g)_j \quad \text{in\,}
U_i\cap U_j.
\end{equation*}

There exist smooth sections $\psi_i$ over $U_i$, vanishing on $D$,
such that
\begin{equation}\label{deltalogg}
(\delta\log g)_{ij}= \psi_i -\psi_j \quad \text{in\,} U_i\cap U_j.
\end{equation}
One may for example choose a smooth function $\rho: M\to [0,1]$
which vanishes in a neighbourhood of $\supp D \cup \supp\gamma_g$
and equals one on each $U_i\cap U_j$, $i\ne j$ and define
$\psi_i = \rho (\log g)_i$ in $U_i.$ In any case, \eqref{deltalogg} shows that the $\psi_i$ satisfy
\begin{equation*}
\bar{\partial}\psi_i = \bar{\partial}\psi_j \quad \text{in\,}
U_i\cap U_j,
\end{equation*}
so that $\{\bar{\partial}\psi_i\}$ defines a global $(0,1)$-form
 $\bar{\partial}\psi$ on $M$.
The Serre pairing is then defined by
\begin{equation*}
\scal{\frac{df}{f}}{[\delta\log g]}_{\mathrm{Serre}} =\frac{1}{2\pi\I}\int_M \frac{df}{f}
\wedge \bar{\partial}\psi.
\end{equation*}
It is straightforward to check that the result ($\Mod 2\pi \I$) does
not depend upon the choices made, and that it ($\Mod 2\pi\I$ )
agrees with $\int_M \frac{df}{f}\wedge d\,\Log g$.

A variant of the above is to consider the product
$\frac{df}{f}\wedge [\delta\log g]$ directly as an element in $H^1
(M, \mathcal{O}^{1,0})$, because there is a natural multiplication map
$$
H^0(M, \mathcal{O}_{D}^{1,0}) \times H^1(M, \mathcal{O}_{-D})\to H^1
(M, \mathcal{O}^{1,0}),
$$
and use the residue map (sum of residues;
see \cite{Gunning72}, \cite{Forster})
$$
\res : H^1 (M, \mathcal{O}^{1,0}) \to \C.
$$
Then one verifies that
$$
\res \,( \frac{df}{f}\wedge [\delta\log g]) =\frac{1}{2\pi\I}\int_M
\frac{df}{f}\wedge d\,\Log g \quad (\Mod 2\pi\I).
$$

In summary we have
\begin{theorem} For any two meromorphic functions $f$ and $g$
$$
\Res(f,g)= \exp(\scal{\frac{df}{f}}{ [\delta\log g]}_{\mathrm{Serre}}) =\exp(\res \,( \frac{df}{f}\wedge [\delta\log g])).
$$
\end{theorem}

The above expressions can be viewed as polarized and global versions of the torsor, or local symbol, as studied by P.~Deligne, see in particular Example~2.8 in \cite{Deligne}.

\section{Determinantal formulas}
\label{sec:seg}

\subsection{Resultant via Szeg\"o's strong limit theorem}
\label{seq:seg}
In this section we show that the resultant of two rational functions on $\mathbb{P}$ admits several equivalent representations, among others as a Cauchy determinant and as a determinant of a truncated Toeplitz operator.
We start with establishing a connection between resultants and  Szeg\"o's strong limit theorem.

Let us apply the results of the previous section to the case when
\begin{equation*}
M={\mathbb P}, \quad U=\D, \quad V={\mathbb P}\setminus\overline{\D},\quad
\Gamma= {\mathbb T}\equiv \partial\D,
\end{equation*}
and $h$ is holomorphic and nonvanishing
in a neighbourhood of ${\mathbb T}$ with
$\wind_\mathbb{T}h=0$ (equivalent to that $\log h$ has a single-valued branch on  ${\mathbb T}$ in this case). Choose an arbitrary branch, $\Log\,h$, and expand it  in a Laurent series
\begin{equation*}
\Log\,h (z)= \sum_{-\infty}^\infty s_k z^k.
\end{equation*}
Note that $s_0$ is determined modulo $2\pi\I \mathbb{Z}$ only and that the $s_k$
also are the Fourier coefficients of $\Log h(e^{\I \theta})$:
\begin{equation}\label{defa}
s_k=(\Log  h)_k=\frac{1}{2\pi }\int_{0}^{2\pi}e^{-\I k\theta}\Log h(e^{\I \theta})\;d\theta.
\end{equation}

Then using the Cauchy kernel \eqref{Cauchykernel} with $z_0=\infty$
one gets
\begin{equation*}
(\Log\,h)_+ (z)= \sum_{k=0}^\infty s_k z^k,\quad
(\Log\,h)_- (z)= -\sum_{k=1}^{\infty}s_{-k} z^{-k},
\end{equation*}
and $d(\Log\,h)_-(z) =\sum_{k=1}^{\infty} ks_{-k}\frac{dz}{z^{k+1}}.$
This gives the formula
\begin{equation}
\Res_{\mathbb{D}}(h)=\exp[\sum_{k=1}^{\infty} ks_k s_{-k}].
\end{equation}
In particular, we have the following corollary of Theorem~\ref{the:coh}.

\begin{corollary}
\label{pro:cohomol}
Let $f$ and $g$ be two rational functions with
$\supp (f)\subset \mathbb{D}$ and $\supp (g)\subset \mathbb{P}\setminus \overline{\mathbb{D}}$. Then
\begin{equation}
\label{eq:coh}
\Res (f,g)=\Res_\mathbb{D} (\frac{f}{g})=\exp[\sum_{k=1}^{\infty} ks_k s_{-k}],
\end{equation}
where $\Log \frac{f(e^{\I \theta})}{g(e^{\I \theta})}=\sum_{k=-\infty}^\infty s_k e^{\I k \theta}$
is the corresponding  Fourier series.
\end{corollary}

The right member in (\ref{eq:coh}) admits a clear interpretation in terms of the celebrated Szeg\"o strong limit theorem (see \cite{Botch99} and the references therein). Indeed, under the assumptions of Corollary~\ref{pro:cohomol},
\begin{equation*}
h(e^{\I \theta})=\frac{f(e^{\I \theta})}{g(e^{\I \theta})}=\sum_{k=-\infty}^\infty h_k e^{\I k \theta}\in L^\infty(\mathbb{T}),
\end{equation*}
therefore $h$ naturally generates a Toeplitz operator on the Hardy space $H^2(\mathbb{D})$:
\begin{equation*}
T(h): \phi\to\mathrm{ P}_+(h\phi),
\end{equation*}
where $\phi\in H^2(\mathbb{D})$ and $\mathrm{P}_+: L^2(\mathbb{T})\to H^2(\mathbb{D})$ is the orthogonal projection. Denote by $t(h)$ the corresponding (infinite) Toeplitz matirx
\begin{equation*}
t(h)_{ij}=h_{i-j}, \qquad i,j\geq 1
\end{equation*}
in the orthonormal basis $\{e^{\I k \theta}\}_{k\geq 0}$.

Then the Szeg\"o strong limit theorem says  that, after an appropriate normalization, the determinants of truncated Toeplitz matrices
$\det t_N(h)$ (defined by (\ref{days}))
approach a nonzero limit provided $h$ is sufficiently smooth, has no zeros on $\mathbb{T}$ and the winding number vanishes: $\wind_{\mathbb{T}}(h)=0$
 (see \cite{Botch99}, \cite{Sim}).

To be more specific, under the assumptions made, the operator $T(1/h)T(h)$ is of determinant class (see for the definition  \cite[p.~49]{Sim}) and
\begin{equation}\label{one}
\begin{split}
\lim_{N\to \infty}e^{-N(\Log h)_0}\det t_N(h)&=\exp \sum_{k=1}^\infty k(\Log  h)_k(\Log  h)_{-k}
=\det T(1/h)T(h),
\end{split}
\end{equation}
where $(\Log  h)_k=s_k$ are defined by (\ref{defa}).
Thus $\Res_{\mathbb{D}}(h)=\det T(1/h)T(h)$.

We have the following determinantal characterization of the resultant (cf. (\ref{pincus11})).

\begin{proposition}\label{spirit}
Under assumptions of Corollary~\ref{pro:cohomol}, the multiplicative commutator
\begin{equation*}
T(g)T(f)^{-1}T(g)^{-1}T(f)
\end{equation*}

is
of determinant class and
\begin{equation}\label{both}
\begin{split}
\Res(f,g)&=\det T(\frac{f}{g})T(\frac{g}{f})
=\det [ T(f)^{-1}T(g)T(f)T(g)^{-1}]\\
&=\lim_{N\to \infty}\left(\frac{ g(0)}{f(\infty)}\right)^N\cdot \det t_N(\frac{f}{g})\\
&=\exp \sum_{k=1}^\infty k(\Log  h)_k(\Log  h)_{-k}.\\
\end{split}
\end{equation}
\end{proposition}

\begin{proof}
In view of Corollary~\ref{pro:cohomol}, it suffices  only to establish that the operator determinants and the limit in  (\ref{both}) are equal. Assume that $f$ and $g$ are given by (\ref{fgdefnew}). Then
\begin{equation*}
h(z)=\frac{f(z)}{g(z)}=\frac{f(\infty)}{ g(0)}\cdot\prod_{i=1}^{m}\frac{1-\frac{a_i}{z}}{1-\frac{b_i}{z}}\prod_{j=1}^{n}\frac{1-\frac{z}{d_i}}{1-\frac{z}{c_i}}.
\end{equation*}
Expanding the logarithm
\begin{equation*}\label{logh}
\begin{split}
\Log h(z)&=\Log \frac{f(\infty)}{g(0)}+\sum_{i=1}^m \Log\frac{1-a_i/z}{1-b_i/z}+\sum_{j=1}^n \Log\frac{1-z/d_j}{1-z/c_j}
\end{split}
\end{equation*}
in the Laurent series on unit circle $|z|=1$ we obtain: $(\Log h)_0=\Log \frac{f(\infty)}{g(0)}$ and
\begin{equation*}\label{podst}
(\Log h)_k=\frac{1}{k}\cdot
\left\{
\begin{array}{cc}
\sum_{i=1}^m (a_i^{-k}-b_i^{-k}),& \text{if} \quad k<0 \\
\\
\sum_{j=1}^n (c_i^{-k}-d_i^{-k}) & \text{if} \quad k>0 .
\end{array}
\right.
\end{equation*}
By the assumptions on the zeros and poles of $f$ and $g$, this yields that $\sum_{k\in \mathbb{Z}}|k|\cdot|(\Log h)_k|^2<\infty$. By the Widom theorem \cite{Wid} (see also \cite[p.~336]{Sim}) we conclude that $T(h)^{-1}T(h)-I$ is of trace class.
Therefore the Szeg\"o theorem becomes applicable for $h(z)$.
Inserting the found value $(\Log h)_0$ into (\ref{one}) we obtain
\begin{equation*}\label{one1}
\lim_{N\to \infty}\left(\frac{g(0)}{f(\infty)}\right)^N\cdot \det t_N(h)=\det T(1/h)T(h).
\end{equation*}

It remains only to show that
\begin{equation*}
T(1/h)T(h)=T(f)^{-1}T(g)T(f)T(g)^{-1}.
\end{equation*}
In order to prove this, notice that  by our assumptions
$g, 1/g\in H^2(\mathbb{D})$ with $\sup_{z\in\mathbb{D}}\left|g(z)\right|<\infty$,
and
$f(1/z)\in H^2(\mathbb{D})$ with $\inf_{z\in {\mathbb{D}}}|f(1/z)|>0$.  Thus
$h(z)=f(z)/g(z)$ is the Wiener-Hopf factorization (see, for example, \cite{Sim}, Corollary~6.2.3), therefore
$
T(h)=T(f)T(1/g)=T(f)T(g)^{-1}.
$
Similarly we get $T(1/h)=T(f)^{-1}T(g)$ and desired identity follows.
\qed\end{proof}

\subsection{Cauchy identity}\label{sec:Cauchy}
A related expression for the resultant for two rational functions is given in terms of classical Schur polynomials. Namely, the well-known Cauchy identity \cite[p.~299, p.~323]{Stanley} reads as follows:
\begin{equation}\label{CauchyL}
\prod_{i=1}^m\prod_{j=1}^n\frac{1}{1-a_i c_j}=\sum_{\lambda}S_\lambda(a)S_\lambda( c)=\exp \sum_{k= 1}^\infty kp_k(a)p_k( c).
\end{equation}
Here $\lambda=(\lambda_1,\lambda_2,\ldots,\lambda_k,\ldots)$ denotes a partition, that is a sequence of non-negative numbers in decreasing order  $\lambda_1\geq \lambda_2 \geq\ldots$ with a finite sum,
$$
S_\lambda(x)\equiv s_\lambda(x_1,x_2,\ldots)=\frac{\det (x_i^{\lambda_j+m-j})_{1\leq i,j\leq m}}{\det (x_i^{j})_{1\leq i,j\leq m}}=
\frac{\det (x_i^{\lambda_j+m-j})_{1\leq i,j\leq m}}{\prod\limits_{1\leq i<j\leq m}(x_i-x_j)}
$$
stands for the Schur symmetric polynomials and
$$
p_k(a)=\frac{1}{k}\sum_{i=1}^m a_i^k, \quad p_k( c)=\frac{1}{k}\sum_{j=1}^n  {c}_j^k,
$$
are the so-called power sum symmetric functions.

Note that the series in (\ref{CauchyL}) should be understood in the sense of formal series or the inverse limit (see \cite[p.~18]{MacD}). But if we suppose that
\begin{equation}\label{ab}
|a_i|<1,\quad  |c_j|<1,\quad \forall i,j, 
\end{equation}
then the above identities are valid in the usual sense.

Let us assume that (\ref{ab}) holds. In order to interpret (\ref{CauchyL}) in terms of the meromorphic resultant, we introduce two rational functions
$$
f(z)=\prod_{i=1}^m(1-\frac{a_i}{z}), \qquad g(z)=\prod_{j=1}^n(1-zc_i).
$$
We find
\begin{equation*}
\Res(f,g)=\frac{\prod_{i=1}^m g(a_i)}{g(0)^m}=\prod_{i=1}^{m}\prod_{j=1}^{n}(1-a_i c_j),
\end{equation*}
and by comparing with (\ref{CauchyL}) we obtain
\begin{equation}\label{lastt}
\Res(f,g)=\exp[ -\sum_{k=1}^\infty kp_k(a)p_k(c)].
\end{equation}

By virtue of assumption (\ref{ab}), $\supp (f)\in \mathbb{D}$ and $\supp (g)\in \mathbb{P}\setminus \overline{\mathbb{D}}$, which is consistent with Corollary~\ref{pro:cohomol}. One can easily see that (\ref{lastt}) is a particular case of (\ref{eq:coh}).


\section{Application to the exponential transform of quadrature domains}
\label{seq:appl1}

\subsection{Quadrature domains and the exponential transform}
A bounded domain $\Omega$ in the complex plane
is called a (classical) {\it quadrature domain} \cite{Aharonov-Shapiro76},
\cite{Sakai82}, \cite{Shapiro92}, \cite{Gustafsson-Shapiro05}
or, in a different terminology, an {\it algebraic domain} \cite{Varchenko-Etingof94},
if there exist finitely many points $z_i\in \Omega$ and coefficients $c_i\in \C$
($i=1,\dots, N$, say) such that
\begin{equation}\label{QD}
\int_\Omega h\, dxdy =\sum_{i=1}^N c_i h(z_i)
\end{equation}
for every integrable analytic function $h$ in $\Omega$.
(Repeated points $z_i$ are allowed and should be interpreted as the
occurrence of corresponding derivatives of $h$ in the right member.)

An equivalent characterization is due to  Aharonov and Shapiro
\cite{Aharonov-Shapiro76} and (under simplifying assumptions) Davis
\cite{Davis74}: $\Omega$ is a quadrature domain if and only if there
exists a meromorphic function $S(z)$ in $\Omega$ (the poles are
located at the quadrature nodes $z_i$) such that
\begin{equation}\label{Szz}
S(z)=\bar{z} \quad \text{for \,} z\in\partial\Omega.
\end{equation}
Thus $S(z)$ is the {\it Schwarz function} of $\partial\Omega$ \cite{Davis74},
\cite{Shapiro92}, which in the above case is meromorphic in all of $\Omega$.

Now let $\Omega$ be an arbitrary  bounded open set in the complex
plane. The moments of $\Omega$ are the complex numbers:
\begin{equation*}\label{a-def}
a_{mn}=\int_{\Omega}z^m\bar z^n \;dxdy.
\end{equation*}
Recoding this sequence (on the level of formal series) into a new sequence $b_{mn}$ by the rule
\begin{equation*}
\sum_{m,n=0}^\infty \frac{b_{mn}}{z^{m+1}\bar w^{n+1}}
=1-\exp(-\sum_{m,n=0}^\infty \frac{a_{mn}}{z^{m+1}\bar w^{n+1}}), \qquad |z|,|w|\gg 1,
\end{equation*}
reveals an established notion of {\it exponential transform}
\cite{Carey-Pincus74}, \cite{Putinar98},
\cite{Gustafsson-Putinar98}. More precisely, this is the function of
two complex variables defined by
\begin{equation*}
E_\Omega (z,w) =\exp[\frac{1}{2\pi\I}\int_\Omega
\frac{d\zeta}{\zeta -z}\wedge \frac{d\bar{\zeta}}{\bar{\zeta} -\bar {w}}].
\end{equation*}
It is in principle defined in all $\Com{2}$, but we shall discuss it only in $(\Com{}\setminus \overline{\Omega})^2$, where it is analytic/antianalytic.

For large enough $z$ and $w$ we have
\begin{equation*}
E_\Omega (z,w)=
1-\sum_{m,n=0}^\infty \frac{b_{mn}}{z^{m+1}\bar w^{n+1}}.
\end{equation*}

\begin{remark}\label{rem:princ}
The exponential transform admits the following operator theoretic
interpretation, due to J.D.~Pincus \cite{Pincus}. Let $T:H\to H$ be a
bounded linear operator in a Hilbert space $H$, with one rank
self-commutator given by
$$
[T^*,T]=T^*T-TT^*=\xi\otimes \xi,
$$ where $\xi\in H$, $\xi\ne 0$. Then
there is a measurable function $g:\Com{}\to[0,1]$
with compact support such that
\begin{equation}\label{pincus1}
\det [T_zT_w^*T_z^{-1}{T_w^*}{}^{-1}]=\exp[\frac{1}{2\pi\I}\int_\C
\frac{g(\zeta)\;d\zeta \wedge d\bar{\zeta}}{(\zeta-z)(\bar{\zeta}
-\bar {w})}],
\end{equation}
where $T_u=T-uI$. The function $g$ is called the \textit{principal function} of
$T$. Conversely, for any given function $g$ with values in $[0,1]$
there is an  operator $T$ with one rank self-commutator such
that (\ref{pincus1}) holds.
\end{remark}

Let $\Omega$ be an arbitrary bounded domain. In \cite{Putinar96} M.~Putinar proved that the following conditions are
equivalent:
\begin{itemize}
\item[\textbf{a)}] $\Omega$ is a  quadrature domain;

\item[\textbf{b)}]
$\Omega$ is determined by some \textit{finite} sequence $(a_{mn})_{0\leq m,n\leq N}$;

\item[\textbf{c)}] for some positive integer $N$ there holds
\begin{equation*}\label{b-deg}
\det(b_{mn})_{0\leq m,n\leq N} = 0;
\end{equation*}

\item[\textbf{d)}] the function $E_{\Omega}(z,w)$ is rational for $z, w$ large, of the kind
\begin{equation}\label{e1}
E_\Omega(z,w)=\frac{Q(z, w)}{P(z)\overline{P(w)}},
\end{equation}
where $P$ and $Q$ are polynomials;

\item[\textbf{e)}] there is a bounded linear operator $T$ acting on a
Hilbert space $H$, with spectrum equal to $\overline{\Omega}$, with
rank one self commutator $[T^*,T]=\xi\otimes \xi$ ($\xi\in H$) and
such that the linear span $\bigvee_{k\geq 0}T^{*k}\xi$ is finite
dimensional.
\end{itemize}

When these conditions hold then the minimum possible number $N$ in \textbf{b}) and \textbf{c}), the degree of $P$ in \textbf{d}),
and the dimension of $\bigvee_{k\geq 0}T^{*k}\xi$ in \textbf{e}) all coincide with the order of the quadrature domain, i.e., the number $N$ in (\ref{QD}). For $Q$, see more precisely below.

Note that $E_\Omega$ is Hermitian symmetric: $E_\Omega(w,z)=\overline{E_\Omega(z,{w})}$
and multiplicative: if $\Omega_1$ and $\Omega_2$ are disjoint then
\begin{equation}\label{disjo}
E_{\Omega_1\cup \Omega_2}(z,w)=E_{\Omega_1}(z,w)E_{\Omega_2}(z,w).
\end{equation}

As $|w|\to \infty$ one has
\begin{equation}\label{Cauchy0}
E_\Omega(z,w)=1-\frac{1}{\bar w}K_{\Omega}(z)+ \mathcal{O}(\frac{1}{|w|^2})
\end{equation}
with $z\in\Com{}$ fixed, where $K_{\Omega}(z)=\frac{1}{2\pi \I}\int\limits_{\Omega} \frac{d\zeta\wedge
d\bar\zeta }{\zeta-z}
$
stands for the Cauchy transform of $\Omega$. On the diagonal $w=z$ we have $E_\Omega(z,z)>0$ for $z\in \Com{}\setminus \overline{\Omega}$ and
\begin{equation*}\label{exex}
\lim_{z\to z_0} E_\Omega(z,z)=0\;
\end{equation*}
for almost all $z_0\in \partial \Omega$ (see \cite{Gustafsson-Putinar98} for details). Thus the information of $\partial \Omega$ is explicitly encoded in $E_\Omega$.

It is also worth to mention that $1-E_\Omega(z,w)$ is positive definite as a kernel, which implies that when $\Omega$ is a quadrature domain of order $N$ then $Q(z,w)$ admits the following representation \cite{Gustafsson-Putinar00}:
\begin{equation*}
\begin{split}
Q(z,w)=P(z)\overline{P(w)}-\sum_{k=0}^{N-1}P_k(z)\overline{P_k(w}),
\end{split}
\end{equation*}
where $\deg P_k=k$.

In the simplest case, when $\Omega=\mathbb{D}(0,r)$,
the disk centered at the origin and of radius $r$, the Cauchy transform and the Schwarz function coincide and are equal to
$\frac{r^2}{z}$, and
\begin{equation}\label{et3}
E_{\mathbb{D}(0,r)}(z,w)=1-\frac{r^2}{z \bar w}.
\end{equation}

\subsection{The elimination function on a Schottky double} Let $\Omega$ be a finitely connected
plane domain with analytic boundary or, more generally, a bordered Riemann surface and let
\begin{equation*}
M=\widehat\Omega=\Omega \cup\partial\Omega\cup\widetilde{\Omega}
\end{equation*}
be the Schottky double of $\Omega$, i.e., the compact Riemann surface
obtained by completing $\Omega$ with a backside with the opposite
conformal structure, the two surfaces glued together along
$\partial\Omega$ (see \cite{Farkas-Kra}, for example). On
$\widehat\Omega$ there is a natural anticonformal involution $\phi
:\widehat\Omega\to \widehat\Omega$ exchanging corresponding points
on $\Omega$ and $\widetilde{\Omega}$ and having $\partial\Omega$ as fixed
points.

Let $f$ and $g$ be two meromorphic functions on $\widehat\Omega$.
Then
\begin{equation*}
f^*= \overline{ (f\circ \phi}),\quad g^*= \overline{ (g\circ \phi}).
\end{equation*}
are also meromorphic on $\widehat\Omega$.

\begin{theorem}\label{th:eliminationsymmetric}
With $\Omega$, $\widehat{\Omega}$, $f$, $g$  as above, assume in
addition that $f$ has no poles in $\Omega\cup\partial\Omega$ and
that $g$ has no poles in $\widetilde{\Omega}\cup\partial\Omega$. Then, for large $z$, $w$,
\begin{equation*}
\mathcal{E}_{f,g}(z,\bar{w})= \exp[\frac{1}{2\pi \I}\int_\Omega
\frac{df}{f-z}\wedge
\frac{d\overline{g^*}}{\overline{g^*}-{\overline{w}}}
].
\end{equation*}
In particular,
\begin{equation*}
\mathcal{E}_{f,f^*}(z,\bar{w})= \exp[\frac{1}{2\pi \I}\int_\Omega
\frac{df}{f-z}\wedge
\frac{d \overline{f}}{\overline{f}-\overline{w}}].
\end{equation*}
\end{theorem}

\begin{proof}
For the divisors of $f-z$ and $g-w$ we have, if $z,w$ are large
enough, $\supp(f-z)\subset\widetilde{\Omega}$, $\supp(g-w)\subset{\Omega}$.
Moreover, $\log (g-w)$ has a single-valued branch in
$\widetilde{\Omega}$ (because the image $g(\widetilde{\Omega})$ is contained
in some disk $\D(0,R)$, hence $(g-w)(\widetilde{\Omega})$ is contained
in $\D(-w,R)$, hence $\log (g-w)$ can be chosen single-valued in
$\widetilde{\Omega}$ if $|w|>R$). Using that $g=\overline{g^*}$ on
$\partial\Omega$ we therefore get
\begin{equation*}
\begin{split}
\mathcal{E}_{f,g}(z,\bar{w})
&=\exp[\frac{1}{2\pi \I}\int_{\hat{\Omega}} \frac{df}{f-z}\wedge d\,
\Log (g-\bar{w})]
= \exp[\frac{1}{2\pi \I}\int_{{\Omega}}
\frac{df}{f-z}\wedge d\, \Log (g-\bar{w})]\\
&= \exp[-\frac{1}{2\pi \I}\int_{\partial\Omega}
\frac{df}{f-z}\wedge  \Log (g-\bar{w})]= \exp[-\frac{1}{2\pi \I}\int_{\partial\Omega}
\frac{df}{f-z}\wedge  \Log (\bar{g^*}-\bar{w})]\\
&=\exp[\frac{1}{2\pi\I}\int_\Omega \frac{df}{f -z}\wedge
\frac{d\bar{g^*}}{\bar{g^*}-\bar{w}}].
\end{split}
\end{equation*}
as claimed.
\qed\end{proof}

\subsection{The exponential transform as the meromorphic resultant}\label{uvidim}

Let $S(z)$ be the Schwarz function of a quadrature domain $\Omega$.
Then the relation \eqref{Szz} can be interpreted as saying that the
pair of functions $S(z)$ and $\bar{z}$ on $\Omega$ combines into a
meromorphic function on the Schottky double
$\widehat\Omega=\Omega\cup\partial\Omega\cup\widetilde{\Omega}$ of
$\Omega$, namely
the function $g$ which equals $S(z)$ on $\Omega$, $\bar{z}$ on
$\widetilde{\Omega}$.

The function $f=g^*=\overline{g\circ\phi}$ is then
represented by the opposite pair: $z$ on $\Omega$, $\overline{S(z)}$
on $\widetilde{\Omega}$. It is known \cite{Gustafsson83} that $f$ and
$g=f^*$ generate the field of meromorphic functions on
$\widehat\Omega$, and we call this pair the \textit{canonical representation}
of $\Omega$ in $\widehat{\Omega}$

From Theorem~\ref{th:eliminationsymmetric} we immediately get
\begin{theorem}\label{th:m}
For any quadrature domain $\Omega$
\begin{equation*}
E_\Omega (z,w) = \mathcal{E}_{f,f^*}(z,\bar{w}) \qquad
(|z|,|w|\gg 1),
\end{equation*}
where $f$, $f^*$ is the canonical representation of $\Omega$ in
$\widehat{\Omega}$.
\end{theorem}

Here we used Theorem~\ref{th:eliminationsymmetric} with
$f(\zeta)=\zeta$ on $\Omega$, i.e., $f|_\Omega={\rm id}$. A slightly
more flexible way of formulating the same result is to let $f$ be
defined on an independent surface $W$, so that $f:W\to \Omega$ is a
conformal map. Then $\Omega$ is a quadrature domain if and only if
$f$ extends to a meromorphic function of the Schottky double
$\widehat{W}$ (this is an easy consequence of \eqref{Szz}; cf.
\cite{Gustafsson83}). When this is the case the exponential
transform of $\Omega$ is
\begin{equation*}
E_\Omega (z,w) = \mathcal{E}_{f,f^*}(z,\bar{w}),
\end{equation*}
with the elimination function in the right member now taken in
$\widehat{W}$.

\begin{remark}\label{rem84}
If $\Omega$ is simply connected one may take $W=\D$, so that
$\widehat{W}=\mathbb{P}$ with involution $\phi:\zeta\mapsto
1/\bar{\zeta}$. Then $f:\D\to \Omega$ is a rational function
when (and only when) $\Omega$ is a quadrature domain, hence we
conclude that $E_\Omega (z,w)$ in this case is the elimination
function for two rational functions, $f(\zeta)$ and
$f^*(\zeta)=\overline{f(1/\bar{\zeta}})$. This topic will be
pursued in Section~\ref{seq86}.
\end{remark}

In analogy with (\ref{utvidg}) one can also introduce an extended
version of the exponential transform:
\begin{equation*}\label{utexp}
E_\Omega (z,w;z_0,w_0) :=\exp[\frac{1}{2\pi\I}\int_\Omega
\biggl(\frac{d\zeta}{\zeta -z}-\frac{d\zeta}{\zeta -z_0}\biggr)\wedge \biggl(\frac{d\bar\zeta}{\bar\zeta -\bar  w}-\frac{d\bar\zeta}{\bar\zeta - \bar w_0}\biggr)].
\end{equation*}

One advantage with this \textit{extended exponential transform} is that it is defined for a wider class of domains, for example, for the entire complex plane. If the standard exponential transform is well-defined then
\begin{equation*}\label{quadr-ex}
E_\Omega(z,w;z_0,w_0)=\frac{E_\Omega(z,w)E_\Omega(z_0,w_0)}{E_\Omega(z,w_0)E_\Omega(z_0,w)}.
\end{equation*}
In other direction, the standard exponential transform can be  obtained from the extended version by passing to the limit:
\begin{equation*}
E_\Omega (z,w)=\lim_{z_0,w_0\to\infty}E_\Omega (z,w;z_0,w_0).
\end{equation*}

Arguing as in the proof of Theorem~\ref{th:m} we obtain the following generalization.

\begin{corollary}Let $\Omega$ is a quadrature domain with canonical
representation $f$ and $f^*$. Then
\begin{equation*}\label{th3ext}
E_\Omega (z,w;z_0,w_0)=\mathcal{E}_{f,f^*} (z, \bar w;z_0, \bar w_0),
\end{equation*}
where $\mathcal{E}_{f,f^*} (z,w;z_0,w_0)$ is the extended elimination function (\ref{utvidg}).
\end{corollary}


\subsection{Rational maps}
Now  we study how the exponential transform of an \textit{arbitrary} domain in $M=\mathbb{P}$ behaves under rational maps. For simplicity,  we  only deal with bounded domains, but this restriction is not essential. It can be easily removed by passing to the extended version of the exponential transform.

For domains in general, the exponential transform need not be rational. However we still have the limit relation (\ref{Cauchy0}). This makes it possible to continue $E_\Omega$ at infinity by
\begin{equation*}\label{conti}
E_\Omega(z,\infty)=E_\Omega(\infty,w)=E_\Omega(\infty,\infty)=1.
\end{equation*}

\begin{theorem}\label{TH1}
Let $\Omega_i$, $i=1,2$, be two bounded open sets in the complex plane and $F$ be a $p$-valent proper
rational function which maps  $\Omega_1$ onto $\Omega_2$. Then for all $z,w\in\Com{}\setminus\overline{\Omega}_2$
\begin{equation}\label{www}
\begin{split}
E_{2}^p(z,w)&=E_{1}((F-z),(F-w))=\Res_u (F(u)-z,\Res_v(F(v)-w,E_1(u,v))),
\end{split}
\end{equation}
where $E_k=E_{\Omega_k}$. (See (\ref{hhh}) for the notation.)
\end{theorem}

\begin{proof}
We have
\begin{equation*}
\begin{split}
E_{2}^p(z,w)&=\exp (\frac{p}{2\pi \I}\int\limits_{\Omega_2}
\frac{d\zeta\wedge d\overline{\zeta} }{(\zeta-z)(\overline{\zeta}- \bar {w})})
=\exp \biggl(\frac{1}{2\pi \I }\int\limits_{\Omega_1}
\frac{F'(\zeta)\overline{F'(\zeta)}\;d\zeta\wedge d\bar
\zeta}{(F(\zeta)-z)(\overline {F(\zeta)}- \bar w)}\biggr).
\end{split}
\end{equation*}

Let $D_u$ denote the divisor of $F(\zeta)-u$. Then
\begin{equation*}
\frac{F'(\zeta)}{F(\zeta)-z}=\frac{d}{d\zeta}\log
(F(\zeta)-z)=\sum_{\alpha\in \mathbb{P}}\frac{D_z(\alpha)}{\zeta-\alpha},
\end{equation*}
where the latter sum is finite. Conjugating both sides in this identity for $z=w$ we get
\begin{equation*}
\frac{ \overline{F'(\zeta)}}{\overline{F(\zeta)}- \bar{w}}=
\sum_{\beta\in \mathbb{P}}\frac{D_w(\beta)}{\overline{ \zeta}-\overline{\beta}},
\end{equation*}
therefore,
\begin{equation*}
\begin{split}
\frac{F'(\zeta)\overline{F'(\zeta)}}{(F(\zeta)-z)(\overline
{F(\zeta)}- \bar w)}&= \sum_{\alpha\in \mathbb{P}}\sum_{\beta\in \mathbb{P}}\frac{D_z(\alpha)D_w(\beta)}{(\zeta-\alpha)(\overline{ \zeta}-\overline{\beta})}.
\end{split}
\end{equation*}

By assumptions, $F(\zeta)-u$ is different from $0$ and $\infty$ for any choice of $u\in \Com{}\setminus \overline{\Omega}_2$ and
$\zeta \in {\overline{\Omega}}_1$. Hence $\supp D_u\subset \Com{}\setminus\overline{\Omega}_1$.
Thus successively taking the integral over $\Omega_1$ and the exponential gives
\begin{equation*}\label{rhs}
\begin{split}
E_{2}^p(z,w)=\prod_{\alpha,\beta \in \mathbb{P}}E_1(\alpha,\beta)^{D_z(\alpha) D_w(\beta)}=E_{1}(D_z,D_w),
\end{split}
\end{equation*}
which is the first equality in (\ref{www}). Applying (\ref{subseq}) we get the second equality.
\qed\end{proof}

Since the exponential transform is a hermitian symmetric function of its arguments, a certain care is needed when using formula (\ref{www}). The lemma below shows that the meromorphic resultant is merely Hermitian symmetric when one argument is anti-holomorphic.
Indeed, suppose, for example, that $f$ is holomorphic  and $g$ is anti-holomorphic, that is $g(z)=\overline{h(z})$, where $h$ is a holomorphic function. Note that $(g)=(h)$. Therefore
\begin{equation*}
\begin{split}
\Res(g,f)&=f((g))=f((h))=h((f))=\overline{g((f)})=\overline{\Res(f,g}).
\end{split}
\end{equation*}
In summary we have

\begin{lemma}\label{rem:conj}
Let $f(z)$ be holomorphic (or anti-holomorphic) and $g(z)$ be anti-ho\-lo\-mor\-phic (holomorphic resp.) in $z$. Then
\begin{equation}\label{hermsym}
\Res(g,f)=\overline{\Res(f,g)}.
\end{equation}
\end{lemma}

\begin{corollary}
Under the conditions of Theorem~\ref{TH1}, if $E_1$ is rational then $E_2^p$ is also rational.
\end{corollary}

\begin{proof}
First consider the inner resultant $\Res_v(\cdot,\cdot)$ in (\ref{www}). Since $E_1(u,v)$ and $F(v)-w$ are rational and $E_1$ is hermitian, the resultant is a rational function in $u$ and $\bar w$ by virtue of (\ref{def1}) and Sylvester's representation (\ref{sfor}) (see also Lemma~\ref{rem:conj}). Repeating this for $\Res_u(\cdot, \cdot)$ we get the desired property.
\qed\end{proof}

\begin{remark}
The fact that rationality of the exponential transform
is invariant under the action of rational maps is not
essentially new. In the separable case, that is when $E_{\Omega_1}$ is given by a formula like (\ref{e1}), and in addition
$f$ is a one-to-one mapping, the rationality of $E_{\Omega_2}$ was proven by M.~Putinar
(see Theorem~4.1 in \cite{Putinar96}). This original proof used existence of the principal function  (see Remark~\ref{rem:princ}).
\end{remark}

\subsection{Simply connected quadrature domains}\label{seq86}
Even for quadrature domains, Theorem~\ref{TH1} provides a new effective tool for computing  the exponential transform and, thereby, gives explicit information about the complex moments, the Schwarz function etc.

Suppose that $\Omega$ is a simply connected bounded domain and $F$ is a uniformizing map from the unit disk $\mathbb{D}$ onto $\Omega$. P.~Davis \cite{Davis74} and D.~Aharonov and H.S.~Shapiro \cite{Aharonov-Shapiro76} proved that $\Omega$ is a quadrature domain if and only if $F$ is a rational function. The we have (cf. Remark~\ref{rem84}).

\begin{theorem}\label{cor:mmm}
Let $F$ be a univalent rational map of the unit disk onto a bounded domain $\Omega$. Then
\begin{equation}\label{ww8}
\begin{split}
E_\Omega(z,w)&=\Res_u(F(u)-z,F^*(u)-\bar w)
\end{split}
\end{equation}
where $F^*(u)=\overline{F(\frac{1}{\bar{u}}})$.
\end{theorem}

\begin{proof}
We have from (\ref{et3}) that $E_{\mathbb{D}}(u,v)=1-\frac{1}{u \bar v}$. Hence  $E_{\mathbb{D}}(u,\cdot)$ has a zero at $\frac{1}{\bar u}$ and a pole at the origin, both of order one. Applying (\ref{hermsym}) we find
\begin{equation*}
\begin{split}
\Res_v( F(v)-w,E_{\mathbb{D}}(u,v))&=\overline{\Res_v(E_{\mathbb{D}}(u,v), F(v)-w})
=\frac{\overline{ F(\frac{1}{\bar u}})-\bar w }{\overline{F(0})- \bar w }
=\frac{F^*(u)-\bar w}{\overline{F(0})- \bar w }.
\end{split}
\end{equation*}
Taking into account the null-homogeneity (\ref{homm}) of resultant and using Theorem~\ref{TH1} we obtain (\ref{ww8}).
\qed\end{proof}

Applying (\ref{def1}) can we write the resultant in the right hand side of  (\ref{ww8}) explicitly.

\begin{corollary}\label{remExp}
Let $F(\zeta)=\frac{A(\zeta)}{B(\zeta)}$ be a univalent rational map of the unit disk onto a bounded domain $\Omega$, where $B$ is normalized to be a monic polynomial. Then
\begin{equation}\label{ww80}
\begin{split}
E_\Omega(z,w)=\Res_{\mathrm{pol}}(B,B^\sharp)\cdot \frac{\Res_{\mathrm{pol}}(P_z,P^\sharp_{w})}{T(z)\overline{T(w})},
\end{split}
\end{equation}
where  $m=\deg B$, $n=\max (\deg A,\deg B)=\deg F$, $P_t=A-tB$,
$$
T(z)=(F(0)-z)^{n-m}\Res_{\mathrm{pol}}(P_z,B^\sharp),
$$
and $P^\sharp(\zeta)=\zeta^{\deg P}\overline{P(1/\bar\zeta})$ is the so-called reciprocal polynomial.

\end{corollary}

We finish this section by demonstrating some concrete examples.
First we apply the above results to  polynomial domains. Let, in Corollary~\ref{remExp}, $F(\zeta)=a_1\zeta+\ldots+a_n\zeta^n$ be a polynomial. Then $B=B^\sharp\equiv 1$, $T(z)=z^n$ and
\begin{equation*}
\begin{split}
P_z(\zeta)=-z+a_1\zeta+\ldots+a_n\zeta^n, \qquad
P_w^\sharp(\zeta)&=\bar{a}_n+\ldots+\bar{a}_1\zeta^{n-1}-\bar{w}\zeta^n.
\end{split}
\end{equation*}

This gives the following closed formula.
\begin{equation}\label{Exp1}
\begin{split}
E_{\Omega}(z,w)=\det
\begin{pmatrix}
  -1      &       &       &  \frac{\bar a_n}{\bar w}  &        &\\
  \frac{a_1}{z}   &\ddots &       & \vdots             & \ddots       &\\
  \vdots &       &  -1    &  \frac{\bar a_1}{\bar w}  &  & \frac{\bar a_1}{\bar w}\\
  \frac{a_n}{z}   &       & \frac{a_1}{z}  & -1                 &        &\vdots    \\
         & \ddots& \vdots&                    &\ddots        & \frac{\bar a_1}{\bar w}   \\
         &       & \frac{a_n}{z}  &                    &  & -1\\
\end{pmatrix}.
\end{split}
\end{equation}
A similar determinantal representation is valid also for general rational functions $F$.

For $n=1$ and $n=2$, (\ref{Exp1}) becomes
\begin{equation*}
\begin{split}
E_{\Omega}(z,w)&=1-x_1y_1,\\
E_{\Omega}(z,w)&=1-x_1y_1-2x_2y_2-x_2^2y_2^2-x_1x_2y_1y_2+x_1^2y_2+x_2y_1^2,
\end{split}
\end{equation*}
where $x_i=a_i/z$ and $y_i=\bar a_i/\bar w$.

The determinant in (\ref{Exp1}), and, more generally, the resultant in (\ref{ww8}), has the following transparent interpretation in terms of the Schwarz function. Suppose that $\Omega=F(\mathbb{D})$ for a rational function $F$ and recall the definition (\ref{Szz}) of the Schwarz function of $\partial \Omega$: $S(z)=\bar z$, $z\in \partial \Omega.$
After substitution $z=F(\zeta)$, $|\zeta|=1$, this yields
$$
S(F(\zeta))=\overline{F(\zeta)}=\bar F(\frac{1}{\zeta})=F^*(\zeta).
$$
Note that $F^*(\zeta)$ is a rational function again.
Thus the Schwarz function may be found by elimination of the variable $\zeta$ in the following system of rational equations:
\begin{equation}\label{implicit}
\left\{
\begin{split}
w&=F^*(\zeta),\\
z&= F(\zeta),
\end{split}
\right.
\end{equation}
where $w=S(z)$. Namely, by Proposition~\ref{pr:cr} the system (\ref{implicit}) holds for some $\zeta$ if and only if
\begin{equation}\label{lefthand}
\Res_\zeta(F(\zeta)-z, F^*(\zeta)-w)=0.
\end{equation}

The latter provides an implicit equation for $w=S(z)$ in terms of $z$. Note that the expression on the left hand side in (\ref{lefthand}) is exactly the exponential transform $E_{\Omega}(z,\bar w)$ in (\ref{ww8}). In fact, Theorem~\ref{th:m} implies that for \textit{any}  quadrature domain $\Omega$ one has $E_{\Omega}(z,\overline{S(z)})=0$.

\section{Meromorphic resultant versus polynomial}\label{sec:local}

Recall that the meromorphic resultant vanishes identically  for polynomials (considered as meromorphic functions on  $\mathbb{P}$). This makes it natural to ask whether there is any reasonable reduction of the meromorphic resultant to the polynomial one. Here we shall discuss this question and show how to adapt the main definitions to make them sensible in the polynomial case.

First we recall the concept of  local symbol (see, for example, \cite{Serre59}, \cite{Tate}). Let $f,g$ be meromorphic functions on
an arbitrary Riemann surface $M$. Notice that for any $a\in M$,
the limit
$$
\tau_{a}(f,g):=(-1)^{\ord_a f \ord_a g}\lim_{z\to a} \frac{f(z)^{\ord_a g}}{g(z)^{\ord_a f}}
$$
exists and it is a nonzero complex number. This number
is called the \textit{local symbol} of $f,g$ at $a$.

For all but finitely many $a$ we have $\tau_{a}(f,g)=1$. The following properties follow from the definition:
\begin{equation}\label{symtau}
\tau_{a}(f,g)\tau_{a}(g,f)=1,
\end{equation}
multiplicativity
\begin{equation}\label{multau}
\tau_{a}(f,g)\tau_{a}(f,h)=\tau_{a}(f,gh),
\end{equation}
and
\begin{equation}\label{jac1}
\tau_{a}(f,g)^{\ord_a h}\tau_{a}(g,h)^{\ord_a f}
\tau_{a}(h,f)^{\ord_a g}=(-1)^{{\ord_a f}\cdot{\ord_a g}\cdot{\ord_a h}}.
\end{equation}

In this notation, Weil's reciprocity law in its full strength states that on a \textit{compact} $M$,
the product of the local symbols of any two meromorphic
functions $f$ and $g$ equals one:
\begin{equation}\label{Weil0}
\prod_{a\in M}\tau_{a}(f,g)=1.
\end{equation}

\begin{definition}
\label{def22}
Let $a\in M$ and let $f$ and $g$  be two meromorphic functions which are admissible on $M\setminus \{a\}$. Let $\sigma=\sigma(\zeta)$ be a local coordinate at $a$ normalized such that $\sigma(a)=0$.
Then the following product is well-defined:
\begin{equation}\label{produ}
\Res_{\sigma} (f,g) = \frac{\tau_{a}(\sigma,g)^{\ord_a f}}{\tau_a(f,g)}\prod_{\xi\ne a}g(\xi)^{\ord_\xi f}
\end{equation}
and is called the \textit{reduced} (with respect to $\sigma$) resultant.
\end{definition}

\begin{proposition}
\label{pr:res}
Under the above assumptions,
\begin{equation}\label{ress1}
\Res_{\sigma} (f,g)=(-1)^{\ord_a f \,\ord_a g}\cdot\Res_{\sigma} (g,f) ,
\end{equation}
and
\begin{equation}\label{ress2}
\Res_{\sigma} (f_1f_2,g)=\Res_{\sigma} (f_1,g)\Res_{\sigma} (f_2,g).
\end{equation}
Moreover, if $\sigma'$ is another local coordinate with $\sigma'(a)=0$, then
\begin{equation}\label{ress3}
\Res_{\sigma'} (f,g)=(-\tau_{\xi}(\sigma',\sigma))^{\ord_a f\ord_a g}
\Res_{\sigma} (f,g).
\end{equation}

\end{proposition}

\begin{proof}
Note first $\Res_{\sigma} (f,g)$ vanishes or equals infinity  if and only if $\Res_{\sigma} (g,f)$ does so. Indeed, let us assume that, for instance, $\Res_{\sigma} (f,g)=0$. Then it follows from (\ref{produ}) and the fact that $\tau_{a}(\cdot,\cdot)$ is finite and never vanishes, that $g(\xi_0)^{\ord_{\xi_0}(f)}=0$ for some $\xi_0\ne a$. Hence  $\ord_{\xi_0}(f)\ord_{\xi_0}(g)> 0$, and $f(\xi_0)^{\ord_{\xi_0}(g)}=0$. From the admissibility condition we know that the product $\ord_\xi(f)\ord_\xi(g)$ does not change sign on $M\setminus \{a\}$, therefore $\ord_\xi(f)\ord_\xi(g)\geq 0$ everywhere. Then changing roles of $f$ and $g$ in (\ref{produ}), we get $\Res_{\sigma} (g,f)=0$.

Thus without loss of generality we may assume that
$\Res_{\sigma} (f,g)\ne0$ and $\Res_{\sigma} (f,g)\ne\infty$.
By virtue of the definition of admissibility we see that the product $\ord_\xi f \ord_\xi g$ is semi-definite on $M\setminus\{a\}$, hence
\begin{equation}\label{identic}
\ord_\xi f \ord_\xi g= 0 \qquad (\xi \in M\setminus\{a\}).
\end{equation}

Since $\ord_a \sigma=1 $, we have by (\ref{jac1}) and (\ref{symtau})
$$
\frac{
\tau_{a}(\sigma,f)^{\ord_a g}}{\tau_{a}(\sigma,g)^{\ord_a f}}=
\tau_{a}(g,\sigma)^{\ord_a f}\tau_{a}(\sigma,f)^{\ord_a g}=(-1)^{{\ord_a f}{\ord_a g}}\tau_a(g,f)
$$
We have
\begin{equation*}
\begin{split}
\frac{\Res_{\sigma} (g,f)}{\Res_{\sigma} (f,g)} &=\frac{\tau_a(f,g)
\tau_{a}(\sigma,f)^{\ord_a g}}{\tau_a(g,f)\tau_{a}(\sigma,g)^{\ord_a f}}
\prod_{\xi\ne a}\frac{f(\xi)^{\ord_\xi(g)}}{g(\xi)^{\ord_\xi(f)}}\\
&=(-1)^{{\ord_a f}{\ord_a g}}\tau_a(f,g)\prod_{\xi\ne a}\frac{f(\xi)^{\ord_\xi(g)}}{g(\xi)^{\ord_\xi(f)}}\\
&=(-1)^{{\ord_a f}{\ord_a g}}\tau_a(f,g)\prod_{\xi\ne a}(-1)^{\ord_\xi f \ord_\xi g}\tau_{\xi}(f,g).
\end{split}
\end{equation*}
Hence, by virtue of (\ref{identic}) and (\ref{Weil0}) we obtain
$$
\frac{\Res_{\sigma} (g,f)}{\Res_{\sigma} (f,g)}=(-1)^{{\ord_a f}{\ord_a g}}\prod_{\xi\in M}\tau_{\xi}(f,g)=(-1)^{{\ord_a f}{\ord_a g}},
$$
and (\ref{ress1}) follows.

In order to prove (\ref{ress2}), it suffices to notice that the right side of (\ref{produ}) is multiplicative, by virtue of (\ref{multau}), with respect to $f$.

Finally, we notice that by (\ref{jac1}):
$
\tau_{a}(\sigma',g)\tau_{a}(g,\sigma)
\tau_{a}(\sigma,\sigma')^{\ord_a g}=(-1)^{\ord_a g},
$
hence
\begin{equation*}
\frac{\Res_{\sigma'} (f,g)}{\Res_{\sigma} (f,g)}=\left(\frac{\tau_{a}(\sigma',g)}{\tau_{a}(\sigma\phantom{{}'},g)}\right)^{\ord_a f}=
(-\tau_{a}(\sigma',\sigma))^{\ord_a g\,\ord_a f}
\end{equation*}
and the required formula (\ref{ress3}) follows.
\qed\end{proof}


Now we apply some of the above constructions to the polynomial case. On the Riemann sphere, $\mathbb{P}$, we pick the distinguished point $a=\infty$ and the corresponding local coordinate $\sigma(z)=\frac{1}{z}$.  Since any two polynomials form an admissible pair on
${\mathbb{C}}$, the corresponding product in (\ref{produ}) is well-defined.

Let us consider two arbitrary polynomials  $f$ and $g$.
Since $\ord_\xi f \cdot \ord_\xi g\geq 0$ for any point $\xi$, we see that $\Res_{\sigma} (f,g)=0$ if and only if $f$ and $g$ have
a common zero in ${\mathbb{C}}$. In particular, $\Res_{\sigma}(f,g)\ne 0$ for coprime polynomials.

Now let $f$ and $g$ have no common zeros. In the notation of (\ref{PQ}) we have $\ord_\infty g=-n$ and
\begin{equation*}\label{tauin}
\tau_{\infty}(\sigma, g)=(-1)^{n}\lim_{z\to \infty}\frac{z^{\deg g}}{g(z)}=\frac{(-1)^{n}}{g_n}
\end{equation*}
and
\begin{equation*}
\tau_{\infty}(f,g)=(-1)^{nm}\lim_{z\to \infty} \frac{f(z)^{-n}}{g(z)^{-m}}=(-1)^{nm}\frac{g_n^m}{f_m^n}
\end{equation*}
hence
\begin{equation*}
\Res_{\sigma}(f,g)=f_{m}^{n}\prod_{\xi\ne \infty}g(\xi)^{\ord_\xi(f)}=f_{m}^ng_{n}^m \prod_{i=1}^m\prod_{j=1}^n(a_i-c_j)
\end{equation*}
Thus, comparing this with (\ref{res1}), we recover the classical definition of polynomial resultant. We have therefore
proved the following.

\begin{corollary}\label{cor:res1}
Let $M=\mathbb{P}$ and $\sigma(z)=\frac{1}{z}$ be the standard local coordinate at $\infty$. Then
\begin{equation*}
\Res_{\sigma} (f,g)=\Res_{\mathrm{pol}} (f,g).
\end{equation*}
\end{corollary}

A beautiful  interpretation of the product in the right hand side of (\ref{produ}) as a determinant is given in a recent paper of J.-L.~Brylinski and E.~Previato \cite{BrP}. In particular, the authors show that this product is described as the determinant $\det (f,A/gA)$ of the Koszul double complex for $f$ and $g$ acting on $A=H^0(M\setminus \{a\},\mathcal{O})$.


\end{document}